\documentclass[11pt]{article}
\usepackage{amsfonts}

\usepackage{graphics}
\usepackage{indentfirst}
\usepackage{cite}
\usepackage{latexsym}
\usepackage{amsmath}
\usepackage{amssymb}
\usepackage[dvips]{epsfig}
\usepackage{amscd}

\newcommand{\rs}{\rho_s}
\newcommand{\pr}{P(\rho)  }
\newcommand{\prs}{P(\rho_s)  }
\newcommand{\pms}{ \pr-\prs }
\let\oldsection\section
\renewcommand\section{\setcounter{equation}{0}\oldsection}

\allowdisplaybreaks
\def\pf{\it{Proof.}\rm\quad}

\def\R{\mathbb{R}}\def\pa{\partial}

\newcommand\divg{{\text{div}}}
\newtheorem{defn}{Definition}[section]
\newtheorem{thm}{Theorem}[section]
\newtheorem{pro}{Proposition}[section]
\newtheorem{lem}{Lemma}[section]
\newtheorem{re}{Remark}[section]
\newtheorem{theorem}{Theorem}[section]
\newtheorem{remark}{Remark}[section]

\newtheorem{lemma}[theorem]{Lemma}

\newcommand{\ti}{\tilde}

\renewcommand{\div}{ {\rm div }  }

\newcommand{\na}{\nabla }

\renewcommand{\r}{\mathbb{R}}
\newcommand{\bi}{\bibitem}

\newcommand{\bl}{\begin{lemma}}
\newcommand{\el}{\end{lemma}}
\newcommand{\et}{\end{theorem}}
\newcommand{\ga}{\gamma}

\newcommand{\al}{\alpha}
\newcommand{\de}{\delta}
\newcommand{\ve}{\varepsilon}
\newcommand{\la}{\label}

\newcommand{\p}{p(\rho)  }

\newcommand{\bn}{\begin{eqnarray}}
\newcommand{\en}{\end{eqnarray}}
\newcommand{\bnn}{\begin{eqnarray*}}
\newcommand{\enn}{\end{eqnarray*}}

\newcommand{\bnnn}{\begin{eqnarray*}}
\newcommand{\ennn}{\end{eqnarray*}}

\newcommand{\ba}{\begin{aligned}}
\newcommand{\ea}{\end{aligned}}
\newcommand{\be}{\begin{equation}}
\newcommand{\ee}{\end{equation}}

\def\p{\partial}
\def\norm[#1]#2{\|#2\|_{#1}}

\newcommand{\no}{\nonumber\\}
\newcommand{\ep}{\varepsilon}
\newcommand{\n}{\rho}
\newcommand{\si}{\sigma}

\def\la{\label}

\def\na{\nabla}

\begin{document}
\title{ On the motion of three-dimensional compressible isentropic flows with large external potential forces and vacuum\thanks{This work was partially
supported by
  NNSFC (Grant Nos.  10801111, 10971215 \& 10971171),
  the Fundamental Research Funds for the Central Universities (Grant No. 2010121006),
  and the Natural Science Foundation of Fujian Province of China (Grant No.
  2010J05011).
      }}
\author{{Jing Li}\\
\small Institute of Applied Mathematics, AMSS \& Hua Loo-Keng Key Laboratory of Mathematics, \\
\small Chinese Academy of Sciences,  Beijing 100190, P. R. China\\
\small E-mail: ajingli@gmail.com\\[2mm]
{Jianwen Zhang}\\
\small School of Mathematical Sciences, Xiamen University, Xiamen
361005, P. R. China\\
\small E-mail: jwzhang@xmu.edu.cn\\[2mm]
{Junning Zhao}\\
\small School of Mathematical Sciences, Xiamen University, Xiamen
361005, P. R. China\\
\small E-mail: jnzhao@xmu.edu.cn}

\date{}

\maketitle \noindent{\bf Abstract.} We study the global existence
and uniqueness of classical solutions to the three-dimensional
compressible isentropic Navier-Stokes equations with vacuum and
external potential forces which could be arbitrarily large provided
the initial data is of small energy and the unique steady state is
strictly away from vacuum. In particular, the solution may have
large oscillations and contain vacuum states. For the case of
discontinuous initial data, we also prove the global existence of
weak solutions. The large-time behavior of the solution is obtained
simultaneously. It is worthwhile mentioning that the compatibility
condition on the initial data and the regularity condition of the
external potential forces in the present paper are much weaker than
those assumed in the existing literature.

\vskip 2mm


\section{Introduction}
The motion of three-dimensional viscous compressible isentropic
flows occupying a domain $\Omega\subset\R^3$ is governed by the
compressible Navier-Stokes equations:
\begin{eqnarray}
\rho_t+\divg (\rho u)&=&0,\label{1.1}
\\
(\rho u)_t+\divg(\rho u\otimes u)+\nabla P(\rho)&=&\mu\Delta
u+(\mu+\lambda)\nabla\divg u+\rho \ti  f,\label{1.2}
\end{eqnarray}
where the unknown functions $\rho\geq0,u =(u^1,u^2,u^3)$ and
$P(\rho)=A\rho^\gamma$ $(A>0,\gamma>1)$ are the fluid density,
velocity and pressure, respectively, and $\ti f=\ti f(x)$ is the
known external force. The viscosity coefficients $\mu$ and $\lambda$
satisfy the physical restrictions $\mu>0$ and $3\lambda+2\mu\geq 0$.

Let $\Omega=\r^3 $ and $\rho_\infty>0$ be a fixed positive constant.
For the external force in the form:
\begin{equation}
\ti f=\na f\quad{\rm with}\quad f=f(x)\label{1.3}
\end{equation}
we look for the solutions, $(\rho(x,t),u(x,t) ),$ to the Cauchy
problem of (\ref{1.1}),  (\ref{1.2}) with  the far field behavior
\begin{equation}
(\rho,u )(x,t)\to (\rho_\infty,0  )  \quad\mbox{ as
}\;\; |x|\to\infty,\;\; t>0,\label{1.4}
\end{equation}
and the initial data
\begin{equation}
(\rho,u)(x,0)=(\rho_0,u_0)(x),\quad x\in\R^3.\label{1.5}
\end{equation}

The equations (\ref{1.1}), (\ref{1.2}) describe the conservation
laws of mass and momentum, respectively. There has been a lot of
literature on the existence and the large-time behavior of solutions
to the compressible Navier-Stokes equations. The one-dimensional
problem has been extensively studied, see \cite{ka1968,AKM1990}  and
the references therein.
 For the multi-dimensional case, the local existence and uniqueness of
classical solutions was proved in \cite{Na1962,se1} in the absence
of vacuum, and recently, the local strong solutions were studied in
\cite{CCK2004,CK2006,CK2003,SS1993} for the case that the initial
density need not be positive and may vanish in open sets. The global
classical solutions were first obtained by Matsumura-Nishida
\cite{MN1980} when the initial data are close to a non-vacuum
equilibrium in Sobolev space $H^3$, see also \cite{MN1983} for the
exterior problem. Later, Hoff \cite{Ho1995-1,Ho1995-2} considered
the weak solutions without vacuum and external forces for the
discontinuous initial data. The global theory for the
multi-dimensional compressible Navier-Stokes equations with ``large
data" is more delicate. Vaigant and Kazhikov \cite{VK1995} studied
the global existence of classical solutions in two-dimensional
periodic domain when the viscosity coefficients depend on density in
a very specific way and the initial density is strictly positive.
One of the most important breakthrough about the global theory of
``large data" is the work of Lions \cite{Li1998} (see also Feireisl
et al. \cite{FNP2001, F2003}), who first proved the global existence
of weak solutions (the so-called ``finite energy weak solutions") to
the initial/inital-boundary value problem of (\ref{1.1}),
(\ref{1.2}) with generally large initial data when the initial
energy is finite and the adiabatic exponent $\gamma$ is suitably
large (i.e. $\gamma>3/2$).  Recently, under the additional
assumptions that the viscosity coefficients $\mu $ and $\lambda$
satisfy $\mu>\max\{ 4\lambda,-\lambda\}$ and that the  far field
density is away from vacuum, Hoff et al. (cf. \cite{Ho3,hs,ht})
obtained a new type of global weak solutions with small energy,
which have extra regularity information compared with those large
weak ones constructed by Lions \cite{Li1998} and Feireisl
\cite{FNP2001, F2003}. More recently, Huang-Li-Xin \cite{HLX2010}
established the global existence and uniqueness of classical
solutions to the Cauchy problem for the isentropic compressible
Navier-Stokes equations in three-dimensional space with smooth
initial data which are of small energy but possibly large
oscillations; in particular, the initial density is allowed to
vanish, even has compact support.

It has been mentioned in many papers (see, e.g. \cite{FP1999,
MY2001,LM2010}) that the large external forces will significantly
affect the dynamic motion of compressible flows and cause some
serious difficulties in the mathematical study. In the following, we
briefly recall some recent progress on the multi-dimensional
compressible Navier-Stokes equations subject to external forces.
Indeed, when both the initial perturbations and the external forces
are sufficiently small, there have been many studies on the global
existence and the large-time behavior of the compressible
Navier-Stokes equations, see, for example,
\cite{Du2007-1,Du2007-2,MN1983,ST2003,Va1983}, and among others. For
large external forces,   Feireisl and Petzeltov${\rm\acute a}$
\cite{FP1999}, Novotny and Stra${\rm\check{s}}$kraba \cite{NS2000}
proved for different boundary conditions that if the adiabatic
exponent $\gamma$ is larger than $3/2$ and there exists a unique
steady state, then the density of any global weak solution converges
to the steady state density in some $L^q$-norm as time goes to
infinity. Under the assumptions that the adiabatic exponent $\gamma$
is close to $1$ and the external forces satisfy some decay
properties in the far field, Matsumura and Yamagata \cite{MY2001}
studied the global existence and large-time behavior of weak
solutions to the Cauchy problem of (\ref{1.1}), (\ref{1.2}) when the
initial perturbations are suitably small in $L^2\cap L^\infty$ for
density (away from vacuum) and in $H^1$ for velocity. Recently,
under the same smallness assumptions on the initial perturbations,
Li and Matsumura \cite{LM2010} succeeded in removing the smallness
condition on $|\gamma-1|$ and the decay assumptions on the external
force, and thus, improved the Matsumura-Yamagata's result
\cite{MY2001}.

It is worth noting that among the papers
\cite{Ho1995-1,Ho1995-2,MY2001,LM2010} mentioned above, the initial
perturbation of density around a given positive state is small in
$L^\infty$, which in particular implies that the density is
uniformly away from vacuum. However, as emphasized in many papers
related to compressible fluid dynamics
\cite{CCK2004,CK2006,CK2003,Xin1998},
the possible presence of vacuum is one of the major difficulties
when the problems of global existence, uniqueness and regularity of
solutions to the compressible Navier-Stokes equations are concerned.

Thus,  our main aim in this paper
is to establish global well-posedness theorem and to study the
large-time behavior of solutions for the compressible Navie-Stokes
equations (\ref{1.1}), (\ref{1.2}) with large external forces and
vacuum states.

Before stating the main results, we explain the notations and
conventions used throughout this paper. For simplicity, we denote
$$
\int fdx=\int_{\r^3}fdx.
$$
For $1\le p\le \infty $ and integer $k\ge 0$, we adopt the following
simplified notations for the standard homogeneous and inhomogeneous
Sobolev spaces:
$$
\left\{
\begin{array}{l}
L^p=L^p(\R^3),\quad
W^{k,p}=W^{k,p}(\R^3),\quad H^k=W^{k,2},\\[2mm]
D^1= \left. \left\{u\in L^6(\R^3) \,\right| \|{\nabla
u}\|_{L^2}<\infty
\right\},\\[2mm]
D^{1,p}=\left.\left\{u\in L^1_{\rm loc}(\R^3)\,\right|\|{\nabla
u}\|_{L^p}<\infty\right\} .
\end{array}\right.
$$

We first study the stationary problem of (\ref{1.1})--(\ref{1.5}).
In view of \cite[Remark 2.1]{LM2010}, it is known that the smooth
steady solution, $(\n_s(x),u_s(x))$, is unique and $ u_s\equiv 0$.
Hence, we infer from (\ref{1.2}) that the steady state density
$\rho_s=\rho_s(x)$ satisfies
\begin{equation}
\nabla P(\rho_s)=\rho_s\nabla f,\quad
\rho_s(x)\to\rho_\infty\quad{\rm as}\quad |x|\to\infty,\label{1.6}
\end{equation}
which implies that $\rho_s $ is uniquely
determined by
$$
\int_{\rho_\infty}^{\rho_s}\rho^{-1}P'(\rho)d\rho=f(x).
$$
In order to avoid the vacuum states in $\rho_s$, we suppose that
\begin{equation}
-\int_0^{\rho_\infty}\frac{P'(\rho)}{\rho}d\rho<\inf_{x\in\R^3}
f(x)\leq \sup_{x\in\R^3}
f(x)<\int_{\rho_\infty}^\infty\frac{P'(\rho)}{\rho}d\rho.\label{1.7}
\end{equation}

Thus, it follows from (\ref{1.6}) and (\ref{1.7}) that
\begin{pro}\label{pro1.1} Assume that $f\in H^2$  satisfy
(\ref{1.7}). Then the stationary problem (\ref{1.6}) has a unique
solution $\rho_s=\rho_s(x)$ satisfying
\begin{equation}
\rho_s-\rho_\infty\in H^2  ,\quad 0<\underline\rho\leq\inf_{x\in
\R^3}\rho_s(x)\leq \sup_{x\in
\R^3}\rho_s(x)\leq\bar\rho<\infty,\label{1.8}
\end{equation}
where $\underline\rho,\bar\rho$ are two positive constants depending
only on $A,\gamma,\rho_\infty, $ $\inf\limits_{x\in\R^3} f(x)$, and
$\sup\limits_{x\in\R^3} f(x)$. Furthermore, if $f\in W^{2,q}$ with
some $q\in (3,6)$, then
\begin{equation}
\|\nabla\rho_s\|_{H^1\cap W^{1,q} }\leq C ,\label{1.9}
\end{equation}
where $C$ is a positive constant depending only on
$A,\gamma,\rho_\infty,\inf\limits_{x\in\R^3} f(x)$, and
$\|f\|_{H^2\cap W^{2,q} }$.
\end{pro}
\begin{remark}\label{re1.1} For $P(\n)=A\n^\ga$ with $A>0,\ga>1,$ it is easily seen from (\ref{1.6}) and  (\ref{1.7}) that
$$
\rho_s(x)=\left(\rho_\infty^{\gamma-1}+\frac{\gamma-1}{A\gamma}f(x)\right)^{1\over\gamma-1},
$$
which, together with $f\in H^2$, implies that (\ref{1.8}) holds
provided
$$
\inf_{x\in\r^3} f(x)>-\frac{A\gamma}{\gamma-1}\rho_\infty^{\gamma-1} .
$$
\end{remark}
\begin{remark}\label{re1.2} To study the large-time behavior of global weak solutions, the authors in
\cite{MY2001,LM2010} technically required that $f\in H^3$ which
particularly yields $\nabla\rho_s\in L^\infty$ and plays a key role
in their analysis. Here we only assume that $f\in H^2$ which is much
weaker than that in \cite{MY2001,LM2010} and is only used to
guarantee (\ref{1.8}). The additional condition $f\in W^{2,q}\supset
H^3$ will be used to derive the high order estimates needed for the
global classical existence.
\end{remark}

With the steady state $(\n_s,u_s)$ at hand, we can define the
initial energy $C_0$ as follows:
\begin{equation} C_0\triangleq
\int_{\R^3}\left(G(\rho_0)+\frac{1}{2}\rho_0|u_0|^2\right)dx,
\label{1.10}
\end{equation}
where $G(\cdot)$ is the potential energy density defined by
\begin{equation}
G(\rho)\triangleq\int_{\rho_s}^\rho\int_{\rho_s}^r
\frac{P'(\xi)}{\xi} d\xi
dr=\rho\int_{\rho_s}^\rho\frac{P(\xi)-P(\rho_s)}{\xi^2}d\xi.\label{1.11}
\end{equation}
Now, our first result concerning the global existence of classical
solution of (\ref{1.1})--(\ref{1.5}) can be formulated as follows.
\begin{thm}\label{thm1.1}For $q\in (3,6)$, let $f\in H^2\cap W^{2,q}$ satisfy
(\ref{1.7}) and $\rho_s=\rho_s(x)$ be the steady state density of
(\ref{1.6}). For given positive numbers $M$ ($M$ may be arbitrarily
large) and $\tilde\rho\geq \bar\rho+1,$ assume that the initial data
$(\rho_0,u_0)$ satisfy
\begin{equation}
\left(\rho_0-\rho_\infty,P(\rho_0)-P(\rho_\infty)\right)\in H^2\cap
W^{2,q},\quad u_0\in H^2,\label{1.12}
\end{equation}
\begin{equation}
0\leq\inf_{x\in\R^3}\rho_0(x)\leq\sup_{x\in\R^3}\rho_0(x)\leq\tilde\rho,\quad
\|\nabla u_0\|_{L^2}^2\leq M,\label{1.13}
\end{equation}
and that the compatibility condition
\begin{equation}
-\mu\Delta u_0-(\lambda+\mu)\nabla{\rm div}u_0 +\nabla
P(\rho_0)=\rho_0^{1/2} g\label{1.14}
\end{equation}
holds for some $g\in L^2$. Then there exists a positive constant
$\varepsilon>0$, depending only on $\mu$, $\lambda$, $A$, $\gamma$,
$\rho_\infty$, $\tilde\rho$, $M$, $\inf\limits_{x\in \r^3}f(x)$ and
$\|f\|_{H^2\cap W^{2,q} }$, such that if the initial energy
satisfies
\begin{equation}
C_0 \leq \varepsilon,\label{1.15}
\end{equation}
the Cauchy problem (\ref{1.1})--(\ref{1.5}) has a unique global
classical solution $(\rho,u)$, defined on $\R^3\times(0,T]$ for any
$0<T<\infty$ and satisfying
\begin{equation}
0\leq\rho(x,t)\leq 2\tilde\rho\quad{for\;\;all}\quad
x\in\R^3,\;t\geq0,\label{1.16}
\end{equation}
and
\begin{equation}
\left\{
\begin{array}{lll}
\left(\rho-\rho_\infty,P(\rho)-P(\rho_\infty)\right)\in C([0,T];H^2\cap W^{2,q}),\\[2mm]
u\in C([0,T];H^2)\cap L^\infty(\tau,T; H^3\cap W^{3,q}),\\[2mm]
u_t\in L^\infty(\tau,T;H^2)\cap H^1(\tau,T;H^1),
\end{array}\right.\label{1.17}
\end{equation}
for any $0<\tau<T<\infty$. Moreover,  one has the following
large-time behavior: \be\la{1.18}
  \lim_{t\rightarrow \infty}\left( \|\n(\cdot,t)-\rs\|_{L^p} +\|\na u(\cdot,t)\|_{L^r } \right)=0,  \ee
   with any    \be\la{eq1}  p\in  (2 ,\infty),\quad r\in [2,6).
     \ee
\end{thm}

\begin{remark}\label{re1.3} It is clear from (\ref{1.17}) that the solution
obtained in Theorem \ref{thm1.1} becomes a classical one for any
positive time (\cite[Remark 1.1]{HuLi}). Moreover,  although the
solution has small energy, its oscillations could be arbitrarily
large and the interior vacuum states are allowed.\end{remark}

\begin{remark}\label{re1.4}Recently, in the absence of external force, Huang-Li-Xin \cite{HLX2010} (see also \cite{CK2006}) proved the
global existence of classical solutions of the Cauchy problem of
(\ref{1.1}), (\ref{1.2}) with smooth non-negative initial density
under the following compatibility conditions
\begin{equation}
-\mu\Delta u_0-(\lambda+\mu)\nabla{\rm div}u_0 +\nabla
P(\rho_0)=\rho_0 g\label{1.19}
\end{equation}
with $g\triangleq g(x)$ satisfying
\begin{equation}
\rho_0^{1/2}g\in L^2, \,\,\na g\in L^2,\label{1.20}
\end{equation}
which play a key role in the analysis of \cite{CK2006,HLX2010}. It
is worthwhile noting that the compatibility conditions (\ref{1.19}),
(\ref{1.20}) are much stronger than that in
(\ref{1.14}).\end{remark}

It is well-known that the discontinuous solutions (namely, weak
solutions) are fundamental and important in both the physical and
mathematical theory. So, our second aim is to study the global weak
solutions (see Definition \ref{def1.1}) of (\ref{1.1})--(\ref{1.5}).
\begin{defn}\label{def1.1} A pair of functions $(\rho,u)$ is called a weak solution of (\ref{1.1})--(\ref{1.5}),
provided that
$$\rho-\n_\infty\in L_{\rm loc}^\infty(0,\infty;L^2\cap
L^\infty),\quad u\in L^\infty_{\rm loc}(0,\infty;H^1),
$$
and that for all test functions $\psi\in
\mathcal{D}(\R^3\times(-\infty,\infty))$ and $j=1,2,3$,
\begin{eqnarray*}
&&\qquad\quad\int\rho_0\psi(\cdot,0)dx+\int_0^\infty\int\left(\rho\psi_t+\rho
u\cdot\nabla\psi\right)dxdt=0,
\\[2mm]
&&\int\rho_0
u_0^j\psi(\cdot,0)dx+\int_0^\infty\int\left(\rho u^j
\psi_t+\rho u^j u\cdot\nabla\psi+P(\rho)\psi_{x_j}\right)dxdt\\
&&\quad=\int_0^\infty\int\left(\mu\nabla
u^j\cdot\nabla\psi+(\mu+\lambda)({\rm div}
u)\psi_{x_j}\right)dxdt-\int_0^\infty\int \rho f_{x_j}\psi dxdt.
\end{eqnarray*}
\end{defn}

The global existence   of weak solutions of (\ref{1.1})--(\ref{1.5})
with discontinuous initial data can be stated as follows.
\begin{thm}\label{thm1.2} Let $f\in H^2 $ satisfy
(\ref{1.7}) and $\rho_s=\rho_s(x)$ be the steady state density of
(\ref{1.6}). For given positive numbers $M$ ($M$ may be arbitrarily
large) and $\tilde\rho\geq \bar\rho+1,$ assume that the initial data
$(\rho_0,u_0)$ satisfy
\begin{equation}
\left\{
\begin{array}{l}
\left(\rho_0-\rho_\infty, P(\rho_0)-P(\rho_\infty)\right)\in L^2\cap
L^\infty,\quad u_0\in H^1,\\[2mm]
0\leq\inf\limits_{x\in\R^3}\rho_0(x)\leq\sup\limits_{x\in\R^3}\rho_0(x)\leq\tilde\rho,\quad\|\nabla
u_0\|_{L^2}^2\leq M.
\end{array}\right. \label{1.21}
\end{equation}
Then there exists a positive constant $\varepsilon>0$, depending
only on $\mu, \lambda, A, \gamma, \rho_\infty, \tilde\rho, M,
\inf\limits_{x\in \r^3}f(x),$ and $\|f\|_{H^2}$, such that if
\begin{equation}
C_0\leq \varepsilon,\label{1.22}
\end{equation}
the Cauchy problem (\ref{1.1})--(\ref{1.5}) has a global weak
solution $(\rho, u)$ on $\R^3\times(0,\infty)$ in the sense of
Definition \ref{def1.1} satisfying
\begin{eqnarray}
&0\leq\rho(x,t)\leq 2\tilde\rho\quad{for\;\;all}\quad
x\in\R^3,\; t\geq0,\label{1.23} \\[2mm]
&\rho-\rho_\infty\in C([0,\infty); L^2),\quad    \rho u\in
C([0,\infty); H^{-1}),\quad \nabla u\in
L^2(0,\infty;L^2)\label{1.24}
\end{eqnarray}
and for any $p\in(2,\infty),$
 \be\la{1.25}
  \lim_{t\rightarrow \infty}\left( \|\n(\cdot,t)-\rs\|_{L^p} +\|  u(\cdot,t)\|_{L^p\cap L^\infty  } \right)=0.  \ee
\end{thm}
\begin{remark}\label{re1.5} Theorem \ref{thm1.2} extends those results in
\cite{Ho1995-1,Ho1995-2,Ho3,MY2001,LM2010} to the case that both the
vacuum states and the large external forces are involved. Moreover,
the regularity condition $f\in H^2$ on the external force is much
weaker than the one $f\in H^3$ which is technically needed in
\cite{MY2001,LM2010}. Indeed, more regularities of the solutions
away from $t=0$ can be obtained (cf. \cite{Ho1995-1}).
\end{remark}

\begin{remark}\label{re1.6}
Similar to \cite{HLX2010}, the  condition  $\|\nabla
u_0\|_{L^2}^2\leq M $ in both (\ref{1.13}) and  (\ref{1.21}) can be
replaced by  $\|u_0\|_{\dot H^\beta} \le M  $ with any $\beta\in
(1/2,1]$.
\end{remark}

We now comment on the analysis of this paper. Note that  the
compatibility condition (\ref{1.14}) is much weaker than those in
\cite{CK2006}, i.e.,  (\ref{1.19}), (\ref{1.20}) (see Remark
\ref{1.4}), we cannot apply the local existence theorem of classical
solution in \cite{CK2006} to the problem considered. Indeed, we
shall split the proof of Theorem \ref{thm1.1} into three steps.
Roughly speaking, we first use the well-known Matsumura-Nishida's
theorem (see Lemma \ref{lem2.5}) to guarantee the local existence of
classical solutions with strictly positive initial density, then
extend the local classical solutions globally in time just under the
condition that the initial energy is suitably small (see Proposition
\ref{pro5.1}), and finally let the lower bound of the initial
density go to zero. So, to this end, we need some global a priori
estimates which are independent of the lower bound of density. It
turns out that the key issue in the proof is to derive both the
time-independent upper bound of density and the time-dependent
higher order estimates of $(\rho,u)$. To do this, we will borrow
some ideas due from \cite{Ho1995-1,Ho1995-2,HLX2010,LM2010}.
However, because of the arbitrariness of external forces, the
presence of vacuum states and the weaker compatibility condition
(\ref{1.14}), some new difficulties arise and the methods therein
cannot be applied directly.

First, similar to that in \cite{Ho1995-1,HLX2010}, we begin our
proof with the careful initial layer analysis. To do this, we
technically need the following modified ``effective viscous flux":
\begin{equation}
F\triangleq \rho_s^{-1}\left[(\lambda+2\mu)\divg
u-(P(\rho)-P(\rho_s))\right],\label{1.26}
\end{equation}
which was introduced by Li-Matsumura \cite{LM2010} and is different
from the ones in \cite{Ho1995-1,Ho1995-2,HLX2010,MY2001}. Basing on
(\ref{1.2}), (\ref{1.26}) and the standard $L^p$-estimate of
elliptic system, one can derive some subtle connections among the
modified ``effective viscous flux", the gradient and material
derivative of the velocity, and the pressure (see Lemma
\ref{lem3.3}), which are important in the entire analysis,
particularly in closing the time-independent energy estimates stated
in Lemma \ref{lem3.2}. To obtain such connections, we need to deal
with some difficulties induced by the large external forces.  This
will be done by making a full use of the mathematical structure of
the steady state density (see (\ref{1.6})) and adopting an idea due
to Huang-Li-Xin \cite{HLX2006} (see (\ref{3.29}) below), which
enable us to control the terms associated with the pressure and the
large external force by the deviation of the density $\rho$ from the
steady state $\rho_s$. Moreover, compared with
\cite{Ho1995-1,HLX2010}, the properties of $\|\nabla F\|_{L^2}$ and
$\|\nabla {\rm curl} u\|_{L^2}$ are proved in a different manner due
to the fact that $f,\rho_s$ only belong to $H^2$, which also makes
the analysis here need to be more careful than that in
\cite{LM2010}.

However, unlike that in \cite{LM2010}, it seems difficult to use directly this modified ``effective
viscous flux" $F$ in (\ref{1.26})  to prove the
uniform upper bound of the density $\rho$ since we only assume that  $f\in H^2$ which implies that  $ \nabla\rho_s\in  L^s $ with any $ s\in [2,6]$. Indeed,
we overcome this difficulty by replacing     $F$ by
the following standard ``effective viscous flux" (see Lemma
\ref{lem3.8}):
\begin{equation}
\tilde F\triangleq (\lambda+2\mu)\divg u-(P(\rho)-P(\rho_s)),
\label{1.27}
\end{equation}
which is in a similar form as the one defined in
\cite{Ho1995-1,Ho1995-2,Li1998,HLX2010}. A key observation is that
for $r,r_1\in(3,6]$, $r_2\in(6,\infty)$ and $1/r_1+1/r_2=1/r$, one
has
\begin{eqnarray}
\|\nabla \tilde F\|_{L^r}&\leq& C\left(\|\rho\dot u\|_{L^r}+\|(\rho-\rho_s)\nabla f\|_{L^r}\right)\nonumber\\
&\leq& C\left(\|\rho\dot u\|_{L^r}+\|\nabla
f\|_{L^{r_1}}\|(\rho-\rho_s)\|_{L^{r_2}}\right),\label{1.28}
\end{eqnarray}
where `` $\dot{}$ " denotes the
material derivative. Estimate \eqref{1.28} shows that we only need that $\rs$ satisfies  $ \nabla\rho_s \in L^2\cap L^6.$    This is the main difference between $F$ and
$\tilde F$ defined respectively in (\ref{1.26}) and (\ref{1.27}). By
virtue of (\ref{1.28}), we can apply the Sobolev embedding
inequality to derive a desired estimate of $\|F\|_{L^\infty}$, and
thus, prove the pointwise boundedness of density by using the
Zlotnik inequlaity (see Lemma \ref{lem2.2}).

There are also some new difficulties lying in the proof of the
higher order time-dependent estimates. Indeed, to achieve the
estimates on the derivatives of the solutions, we first prove the
important estimates on the gradients of the density and velocity by
solving a logarithm Gronwall inequality in a similar manner as that
in \cite{HLX2010-2,HLX2010}. As a result, one can also easily obtain
the $L^2$-estimates for the second-order derivatives of density,
pressure and velocity. However, due to the weaker compatibility
condition in (\ref{1.14}) (cf. (\ref{1.19}), (\ref{1.20})), the
method used in \cite{HLX2010} cannot be applied any more to obtain
further estimates needed for the existence of classical solutions.
In fact, instead of the $L^2$-method, we succeed in obtaining these
classical estimates by deriving some desired $L^q$-estimates
($3<q<6$) on the higher-order time-spatial derivatives of the
density and velocity, basing on some careful initial-layer analysis
(see Lemmas \ref{lem4.4}--\ref{lem4.6}).

The rest of this paper is organized as follows. In Sect. 2, we first
recall some known inequalities and facts which will be frequently
used in our analysis. In Sect. 3, we derive the key a priori
estimates of the weighted estimates on the gradient and the material
derivative of the velocity and the pointwise upper bound of the
density, all of which are independent of $t$ and will be used to
proved the large-time behavior. In Sect. 4, we prove the
time-dependent estimates on the higher-order norms of the solutions,
which are needed for the existence of classical solutions. Finally,
the main results (i.e. Theorems \ref{thm1.1} and \ref{thm1.2}) will
be proved in Sect. 5.

\section{Auxiliary lemmas}
In this section, we list some elementary inequalities and known
facts which will be used frequently later. We begin with the
well-known Gagliardo-Nirenberg-Sobolev-type inequality (see
\cite{LSU1968}).
\begin{lem}\label{lem2.1} Let $f\in D^1$ and $g\in L^q\cap D^{1,r}$ with $q\in(1,\infty)$ and $r\in (3,\infty)$. Then there
exists a positive constant $C$, which may depend on $q$ and $r$,
such that
\begin{equation}
\|f\|_{L^6}\leq C \|\nabla f\|_{L^2} ,\label{2.1}
\end{equation}
\begin{equation}
\|g\|_{L^\infty}\leq C\|g\|_{L^q}^{q(r-3)/(3r+q(r-3))}\|\nabla
g\|_{L^r}^{3r/(3r+q(r-3))}.\label{2.2}
\end{equation}
\end{lem}

The following Zlotnik inequality, whose proof can be found in
\cite{Zl2000}, will be used later to prove the uniform-in-time upper
bound of the density.
\begin{lem}\label{lem2.2} Assume that the function $y\in W^{1,1}(0,T)$
solves the ODE system:
$$
y'=g(y)+b'(t)\quad{on}\quad[0,T],\quad y(0)=y_0,
$$
where $b\in W^{1,1}(0,T)$ and $g\in C(\R)$. If $g(\infty)=-\infty$
and
\begin{equation}
b(t_2)-b(t_1)\leq N_0+N_1(t_2-t_1)\label{2.3}
\end{equation}
for all $0\leq t_1\leq t_2\leq T$ with some positive constants $N_0$
and $N_1$, then one has
\begin{equation}
y(t)\leq \max\{y_0,\xi^*\}+N_0<+\infty\quad {
on}\quad[0,T],\label{2.4}
\end{equation}
where $\xi^*\in\R$ is a constant such that
\begin{equation}
g(\xi)\leq -N_1\quad for\quad \xi\geq\xi^*.\label{2.5}
\end{equation}
\end{lem}

In order to obtain the time-dependent estimates of $\|\nabla
u\|_{L^\infty}$ and $\|\nabla\rho\|_{L^2\cap L^6}$, we need the
following Beale-Kato-Majda-type inequality, the proof of which can
be found in \cite{HLX2010-2}.
\begin{lem}\label{lem2.3}For
$3<q<\infty$, assume that $\nabla u\in L^2\cap D^{1,q}$. Then there
exists a constant $C>0$, depending only on $q$,  such that
\begin{equation}
\|\nabla u\|_{L^\infty}\leq C\left(\|{\rm div}u\|_{L^\infty}+
\|\nabla\times u\|_{L^\infty} \right)\ln\left(e+\|\nabla^2
u\|_{L^q}\right)+C\|\nabla u\|_{L^2} +C.\label{2.6}
\end{equation}
\end{lem}

Since there is no vacuum state in the far field, it is easy to show
that $u\in L^2$ even that the density is only nonnegative. Indeed,
we have
\begin{lem}\label{lem2.4} Let $\rho_s=\rho_s(x)$ be the steady state density as in Proposition \ref{pro1.1}.
Assume that $(\varrho,v)$ satisfies the following conditions:
$$
0\leq\varrho\leq 2\tilde\rho,\quad (\varrho-\rho_s,\varrho^{1/2}
v)\in L^2,\quad \nabla v\in L^2.
$$
Then there exists a positive constant $C$, depending only on
$\underline\rho,\bar\rho$ and $\tilde\rho$, such that
\begin{equation}
\|v\|_{L^2}\leq
C\left(\|\varrho^{1/2}v\|_{L^2}+\|\varrho-\rho_s\|_{L^2}^{2/3}\|\nabla
v\|_{L^2}\right).\label{2.7}
\end{equation}
\end{lem}
\pf Indeed, it is easy to see that
\begin{eqnarray*}
\int|v|^2dx&\leq& C(\underline\rho)\int\rho_s^2 |v|^2dx\leq
C(\underline\rho)\left(\int\varrho^2
|v|^2dx+\int(\varrho-\rho_s)^2|v|^2dx\right)\\
&\leq& C(\underline\rho,\tilde\rho)\int\varrho |v|^2dx+
C(\underline\rho)\left(\int|\varrho-\rho_s|^3dx\right)^{2/3}
\left(\int|v|^6dx\right)^{1/3}\\
&\leq&
C(\underline\rho,\tilde\rho)\|\varrho^{1/2}v\|_{L^2}^2+C(\underline\rho,\bar\rho,\tilde\rho)\|\varrho-\rho_s\|_{L^2}^{4/3}\|\nabla
v\|_{L^2}^2,
\end{eqnarray*}
which immediately proves (\ref{2.7}).\hfill$\square$

We end this section with the following local well-posedness theorem
of classical solution to the problem (\ref{1.1})--(\ref{1.5}) when
the initial density is strictly away from vacuum (see, e.g.
\cite{Na1962}, and especially Matsumura-Nishida \cite[Theorem
5.2]{MN1980}).
\begin{lem}\label{lem2.5} Assume that the initial data
$(\rho_0,u_0)$ satisfies
\begin{equation}
(\rho_0-\rho_\infty, u_0)\in H^3,\quad
 \inf_{x\in\R^3}\rho_0(x)>0.\label{2.8}
\end{equation}
Then there exist a small time $T_0>0$ and a unique classical
solution $(\rho,u)$ to the Cauchy problem (\ref{1.1})--(\ref{1.5})
on $\R^3\times(0,T_0]$ such that
\begin{equation}
\inf_{(x,t)\in\R^3\times [0,T_0]}
\rho(x,t)\geq\frac{1}{2}\inf_{x\in\R^3}\rho_0(x) ,\label{2.9}
\end{equation}
and
\begin{equation}
\left\{
\begin{array}{l}
\rho-\rho_\infty\in C([0,T_0];H^3)\cap C^1([0,T_0];H^2),
\\[2mm]
u\in C([0,T_0];H^3)\cap C^1([0,T_0];H^1)\cap L^2(0,T_0;H^4),
\end{array}\right.\label{2.10}
\end{equation}
where $T_0>0$ may depend on $\inf\limits_{x\in\R^3}\rho_0(x)$.
\end{lem}

\section{Time-independent a priori estimates}
This section is concerned with the time-independent (weighted)
energy estimates and the uniform upper bound of density, which are
essential for the proofs of Theorems \ref{thm1.1} and \ref{thm1.2}.
To do this, we assume that $(\rho,u)$, defined over $(0,T)$ with
some positive $T>0$, is a smooth solution of the Cauchy problem
(\ref{1.1})--(\ref{1.5}). For simplicity, we introduce the following
functionals:
\begin{equation}
 \Phi_1(T)\triangleq\sup_{0\leq t\leq T}\sigma\int|\nabla
u|^2dx+\int_0^T\sigma\int\rho|\dot u|^2dxdt,\label{3.1}
\end{equation}
\begin{equation}
\Phi_2(T)\triangleq\sup_{0\leq t\leq T}\sigma^3\int\rho|\dot
u|^2dx+\int_0^T\sigma^3\int|\nabla \dot u|^2dxdt\label{3.2}
\end{equation}
and
\begin{equation}
\Phi_3(T)\triangleq \sup_{0\leq t\leq T}\int|\nabla
u|^2dx,\label{3.3}
\end{equation}
where $\sigma(t)\triangleq\{1,t\}$, and the symbol `` $\dot{}$ "
denotes the material derivative $\dot v=v_t+u\cdot\nabla v$.

The main purpose of this section is to prove the following key a
priori estimates.
\begin{pro}\label{pro3.1} Assume that the conditions of Theorem \ref{thm1.2} are satisfied. There
exist two positive constants $\tilde\varepsilon$ and $K$, depending
only on $\mu,\lambda,A,\gamma,\underline\rho, \bar\rho,\tilde\rho,
\inf\limits_{x\in\R^3}f(x),\|f\|_{H^2}$ and $M$, such that if
$(\rho,u)$ is a smooth solution of (\ref{1.1})--(\ref{1.5})
satisfying
\begin{equation}
\left\{
\begin{array}{lll}
0\leq\rho(x,t)\leq 2\tilde\rho\quad{\rm for\;\; all}\quad (x,t)\in\R^3\times[0,T],\\[2mm]
\Phi_1(T)+\Phi_2(T)\leq 2 C_0^{1/2}\quad{\rm and}\quad
\Phi_3(\sigma(T))\leq 3K,
\end{array}\right.\label{3.4}
\end{equation}
then one has
\begin{equation}
\left\{
\begin{array}{lll}
0\leq\rho(x,t)\leq \frac{7}{4}\tilde\rho\quad{\rm for\;\; all}\quad (x,t)\in\R^3\times[0,T],\\[2mm]
\Phi_1(T)+\Phi_2(T)\leq  C_0^{1/2}\quad{\rm and}\quad
\Phi_3(\sigma(T))\leq 2K,
\end{array}\right.\label{3.5}
\end{equation}
provided that the initial energy $C_0$ defined in (\ref{1.10})
satisfies
\begin{equation}
C_0\leq\tilde\varepsilon.\label{3.6}
\end{equation}
\end{pro}

\pf Proposition \ref{pro3.1}
  follows directly from the following Lemmas \ref{lem3.5}, \ref{lem3.6}, and \ref{lem3.8} with $K$ and $\tilde\varepsilon$ being the positive
constants as in Lemmas \ref{lem3.5} and \ref{lem3.8},
respectively.\hfill$\square$

\begin{re}\label{re3.1} We assume throughout this section that the
initial data and the external force only satisfy the conditions of
Theorem \ref{thm1.2}, and hence, the uniform-in-time estimates
derived in this section can be used to study the existence and
large-time behavior of global weak solutions as stated in Theorem
\ref{thm1.2}.
\end{re}

  For notational convenience, throughout this section we
denote by $C$ or $C_i$ ($i=1,2,\ldots$) the generic positive
constants which may depend on
$\mu,\lambda,A,\gamma,\underline\rho,\bar\rho,\tilde\rho,
\inf\limits_{x\in\R^3}f(x),\|f\|_{H^2}$ and $M$,  but not on $T$. We
also sometimes write $C(\alpha)$ to emphasize that $C$ relies on
$\alpha$.

We start the proof with the following standard energy estimate.
\begin{lem}\label{lem3.1} Let $(\rho,u)$ be a smooth solution to (\ref{1.1})--(\ref{1.5}) on $\R^3\times(0,T]$. Then,
\begin{equation}
\sup_{0\leq t\leq T}\int
\left(\frac{1}{2}\rho|u|^2+G(\rho)\right)dx+\int_0^T\int
\left(\mu|\nabla u|^2+(\mu+\lambda)({\rm div} u)^2\right)dxdt\leq
C_0,\label{3.7}
\end{equation}
where $G(\rho)$ is the potential energy density  defined in
(\ref{1.11}).
\end{lem}
\pf Thanks to (\ref{1.6}), the momentum equation (\ref{1.2}) can be
written as
$$
\rho u_t+\rho u\cdot\nabla u+\left(\nabla
P(\rho)-\rho\rho_s^{-1}\nabla P(\rho_s)\right)=\mu\Delta
u+(\mu+\lambda)\nabla\divg u,
$$
which, multiplied by $u$ and integrated by parts over
$\R^3\times(0,t)$, yields
\begin{eqnarray}
&&\frac{1}{2}\left.\int \rho |u|^2 dx\right|_0^t +\int_0^t\int
u\cdot\left(\nabla P(\rho)-\rho\rho_s^{-1}\nabla
P(\rho_s)\right) dxds\nonumber\\
&&\qquad+\int_0^t\int \left(\mu|\nabla u|^2+(\mu+\lambda)(\divg
u)^2\right)dxds=0.\label{3.8}
\end{eqnarray}
After integrating by parts, one infers from (\ref{1.1}) that
\begin{eqnarray*}
&&\int_0^t\int u\cdot\left(\nabla P(\rho)-\rho\rho_s^{-1}\nabla
P(\rho_s)\right)
dxds\nonumber\\
&&\quad=\int_0^t\int \rho u\cdot\nabla
\left(\int_{\rho_s}^\rho\frac{P'(\xi)}{\xi}d\xi\right)dxds\\
&&\quad=-\int_0^t\int \divg(\rho u)
\left(\int_{\rho_s}^\rho\frac{P'(\xi)}{\xi}d\xi\right)dxds\nonumber\\
&&\quad=\int_0^t\int \rho_t
\left(\int_{\rho_s}^\rho\frac{P'(\xi)}{\xi}d\xi\right)dxds=\left.\int
G(\rho)dx\right|_0^t,
\end{eqnarray*}
which, inserted into (\ref{3.8}), leads to the desired estimate in
(\ref{3.7}).\hfill$\square$

It is clear that for all $0\leq\rho\leq 2\tilde\rho$ and
$\underline\rho\leq\rho_s\leq \bar\rho$, there are positive
constants $C_1$ and $ C_2$  depending only on
$\underline\rho,\bar\rho,$ and $\tilde\rho$, such that
$$
C_1(\rho-\rho_s)^2\leq G(\rho)\leq C_2(\rho-\rho_s)^2,
$$
so that, it readily follows from (\ref{3.7}) that
\begin{equation}
\sup_{0\leq t\leq T}\|\rho-\rho_s\|_{L^2}^2\leq C C_0.\label{3.9}
\end{equation}

The next lemma is concerned with the temporary weighted
$L^2$-estimates on the gradient and the material derivatives of the
velocity, the proof of which will be concluded in Lemma \ref{lem3.6}
below. The idea of the proof mainly comes from
\cite{Ho1995-1,HLX2010,LM2010}. However, due to the arbitrariness of
external potential forces, the presence of vacuum states, and
especially, the weaker regularity assumption of the external forces
(i.e. $f\in H^2$), the analysis here needs to be more careful.
\begin{lem}\label{lem3.2} Let $(\rho,u)$ with $\rho\in[0,2\tilde\rho]$ be a smooth solution of (\ref{1.1})--(\ref{1.5}) on $\R^3\times(0,T]$.
Then there exists a positive constant $C$, depending on
$\tilde\rho$, such that
\begin{eqnarray}
\Phi_1(T)&\leq& C\left(C_0+\int_0^T\sigma\|\nabla u\|_{L^3}^3
dt\right),\label{3.10} \\
\Phi_2(T)&\leq& C\left(C_0+\Phi_1(T)+\int_0^T\sigma^3\|\nabla
u\|_{L^4}^4 dt\right).\label{3.11}
\end{eqnarray}
\end{lem}
\pf It follows from (\ref{1.2}) and (\ref{1.6}) that
\begin{eqnarray}
 \rho\dot u-\mu\Delta u-(\mu+\lambda)\nabla\divg
u+ \nabla\left(P(\rho)-P(\rho_s)\right)  =(\rho-\rho_s) \nabla
f.\label{3.12}
\end{eqnarray}
For $m\ge 0$, multiplying (\ref{3.12}) by $\sigma^m\dot u$ in $L^2$
gives
\begin{eqnarray}
\sigma^m\int \rho|\dot{u}|^2dx&=& \mu \sigma^m\int\Delta
u\cdot\dot{u} dx+ (\mu+\lambda)\sigma^m
\int\nabla\divg u\cdot\dot{u} dx\nonumber\\
&&-\sigma^m \int
\dot{u}\cdot\nabla\left(P(\rho)-P(\rho_s)\right)dx+\sigma^m \int
(\rho-\rho_s)
\dot u\cdot\nabla f dx\nonumber\\
&\triangleq& \sum_{i=1}^{4}I_i, \label{3.13}
\end{eqnarray}
where the right-hand side can be estimated term by term as follows.
First, by the definition of the material derivative `` $\dot{}$ "
and integration by parts, we easily get that
\begin{eqnarray}
I_1&=&-\frac{\mu}{2}\left(\sigma^m\|\nabla u\|_{L^2}^2\right)_t +
\frac{\mu}{2}m\sigma^{m-1}\sigma' \|\nabla u\|_{L^2}^2
-\mu \sigma^m \int \p_iu^j\p_i(u^k\p_ku^j)dx\nonumber \\
&\leq&  -\frac{\mu }{2}\left(\sigma^m\|\nabla u\|_{L^2}^2\right)_t +
Cm\sigma^{m-1}\sigma'\|\na u\|_{L^2}^2 + C \sigma^m \|\nabla
u\|_{L^3}^3,\label{3.14}
\end{eqnarray}
and similarly,
\begin{equation}
I_2\leq -\frac{\mu+\lambda}{2}\left(\sigma^m \|\divg
u\|_{L^2}^2\right)_t+ Cm\si^{m-1}\si'\|\na u\|_{L^2}^2 + C \sigma^m
\|\nabla u\|_{L^3}^3.\label{3.15}
\end{equation}

Secondly, noticing that (\ref{1.1}) implies
\begin{equation}
P(\rho)_t+\divg (P(\rho)u)+(\rho P'(\rho)-P(\rho))\divg
u=0,\label{3.16}
\end{equation}
so that,  using (\ref{1.8}), (\ref{2.1}) and (\ref{3.9}), we obtain
after integrating by parts that
\begin{eqnarray}
I_3 &=& \int \left[ \sigma^m (\divg
u)_t\left(P(\rho)-P(\rho_s)\right)
- \sigma^m (u\cdot\nabla u)\cdot\nabla \left(P(\rho)-P(\rho_s)\right)\right]dx\nonumber \\
&= & \frac{d}{dt}\int\sigma^m \text{div}u\left(P(\rho)-P(\rho_s)\right)dx
- m \sigma^{m-1}\si'\int \divg u \left(P(\rho)-P(\rho_s)\right) dx \nonumber\\
&& + \int\sigma^m \left( (\rho P'(\rho)
- P(\rho))(\divg u)^2+   P(\rho)\p_iu^j\p_ju^i + u\cdot \na u\cdot\na \prs \right)dx \nonumber \\
&\le & \frac{d}{dt}\int\sigma^m \text{div}u(\pms)dx
+ m\si^{m-1}\si'\|P(\rho)-P(\rho_s)\|_{L^2}\|\na u\|_{L^2}\nonumber\\
&& +C \si^m\|\na u\|_{L^2}^2 +C\si^m\|u\|_{L^6}  \|\na u\|_{L^2}\|\na\rs\|_{L^3}  \nonumber\\
&\le & \frac{d}{dt}\int\sigma^m \text{div}u(\pms)dx  +C\|\na
u\|_{L^2}^2 +C m^2\si^{2(m-1)}\si'C_0.\label{3.17}
\end{eqnarray}
Analogously, using the fact that $
\rho_t+\divg((\rho-\rho_s)u)+\divg(\rho_s u)=0 $ due to (\ref{1.1}),
we have from (\ref{1.8}), (\ref{2.1}) and (\ref{3.9}) that
\begin{eqnarray}
I_4 &=& \frac{d}{dt}\int\si^m (\rho-\rho_s)  u\cdot\nabla f
dx-m\si^{m-1}\si'\int (\rho-\rho_s)  u\cdot\nabla f dx \nonumber \\
&&+\si^m\int \left[\div(\rho_s u)u\cdot\nabla f-(\rho-\rho_s)u\cdot\na(  u\cdot\nabla f)+(\rho-\rho_s) u\cdot\nabla u \cdot\nabla f\right]dx \nonumber\\
&\le & \frac{d}{dt}\int\si^m (\rho-\rho_s)  u\cdot\nabla f
dx +C m\si^{m-1}\si'\|\rho-\rho_s\|_{L^2} \|u\|_{L^6}\|\nabla f\|_{L^3}\nonumber\\
&&  +C\si^m\int  \left(|u| |\na u||\na f|+|\na\rho_s||u|^2|\na f|+|\rho-\rho_s||u|^2|\na^2f|\right) dx  \nonumber\\
&\le &  \frac{d}{dt}\int\si^m (\rho-\rho_s)  u\cdot\nabla f dx  +C
m\si^{m-1}\si'C_0^{1/2} \|\na u\|_{L^2}+C\si^m
\|u\|_{L^6} \|\na u\|_{L^2} \|\na f \|_{L^3}\nonumber\\
&& +C\sigma^m\|\na \rho_s \|_{L^3}\|u\|_{L^6}^2\|\na f
\|_{L^3}+C\si^m\|\rho-\rho_s\|_{L^6}\|u\|_{L^6}^2\|\na^2 f
\|_{L^2}  \nonumber\\
&\le & \frac{d}{dt}\int\si^m (\rho-\rho_s)  u\cdot\nabla f dx+C\|\na
u\|_{L^2}^2 +C m^2\si^{2(m-1)}\si'C_0.\label{3.18}
\end{eqnarray}

Combining (\ref{3.13})--(\ref{3.15}), (\ref{3.17}) and (\ref{3.18})
leads to \begin{eqnarray} &&\left(\frac{\mu }{2}\si^m\|\nabla
u\|_{L^2}^2+\frac{(\mu+\lambda)}{2}\si^m\|\divg
u\|_{L^2}^2\right)_t+\int \sigma^m
\rho|\dot{u}|^2dx\nonumber\\
&&\quad\le \frac{d}{dt}\int  \sigma^m\left(\divg u\left(P(\rho)-P(\rho_s)\right)+ (\rho-\rho_s)  u\cdot\nabla f\right)dx\nonumber\\
&&\qquad +C m^2\si^{2(m-1)}\si'C_0 + C(1+m\si^{m-1}\si')\|\na
u\|_{L^2}^2+ C \sigma^m \|\nabla u\|_{L^3}^3,\label{3.19}
\end{eqnarray}
where the first term on the right-hand side can be estimated as
follows:
\begin{eqnarray*}
&&\left|\int  \sigma^m\left(\text{div}u(\pms)+ (\n-\rs)
u\cdot\nabla
f\right)dx\right|\\
&&\quad\leq C\sigma^m\left(\|\rho-\rho_s\|_{L^2}\|\nabla
u\|_{L^2}+\|\rho-\rho_s\|_{L^2}\|u\|_{L^6}\|\nabla
f\|_{L^3}\right)\\
&&\quad\leq \frac{\mu}{4}\si^m\|\nabla u\|_{L^2}^2+C\sigma^m C_0.
\end{eqnarray*}
Thus, choosing $m=1$ in (\ref{3.19}), integrating it over $(0,T)$
and using (\ref{3.7}) one gets (\ref{3.10}).

To prove (\ref{3.11}), operating $\sigma^m\dot u^j
\left(\partial_t+\divg(u\cdot)\right)$ to both sides of (\ref{3.12})
and integrating the resulting equations over $\R^3$, we obtain after
summing them up that
\begin{eqnarray}
\sum_{i=1}^3L_i&\triangleq&\sigma^m\int \dot
u^j\left[\partial_t(\rho \dot u^j)+\partial_k(\rho u^k\dot
u^j)\right]dx-\mu \sigma^m\int \dot u^j\left[\Delta
u^j_t+\partial_k(u^k\Delta u^j)\right]dx\nonumber\\
&&-(\mu+\lambda) \sigma^m\int \dot u^j\left[\partial_j\divg
u_{t}+\partial_k(u^k\partial_j\divg
u)\right]dx\nonumber\\
&=&- \sigma^m\int \dot u^j\left[\left(\partial_j
P(\rho)\right)_t+\partial_k\left(u^k\partial_j\left(P(\rho)-P(\rho_s)\right)\right)\right]dx\nonumber\\
&&+ \sigma^m\int \dot
u^j\left[\rho_t\partial_jf+\partial_k\left(u^k(\rho-\rho_s)
\partial_jf\right)\right]dx\triangleq \sum_{i=1}^2R_i.\label{3.20}
\end{eqnarray}

We now estimate each term in (\ref{3.20}). First, by (\ref{1.1}) one
easily gets that
\begin{equation}
L_1=\frac{1}{2}\frac{d}{dt}\int \sigma^m\rho |\dot
u|^2dx-\frac{m}{2}\sigma^{m-1}\sigma'\int \rho |\dot u|^2
dx.\label{3.21}
\end{equation}
Recalling the definition of `` $\dot{}$ ", we deduce from
integration by parts that
\begin{eqnarray}
L_2 &=&\mu\sigma^m\int \left(|\nabla\dot u|^2 -\partial_k\dot
u^j\partial_k\left(u\cdot\nabla u^j\right)+\partial_k\dot
u^ju^k\Delta u^j\right) dx\nonumber\\
&=&\mu \sigma^m\int \left(|\nabla\dot u|^2 -
\partial_k\dot
u^j\partial_ku^l\partial_l u^j-\partial_k\dot
u^ju^l\partial_{kl}^2u^j-\partial_{kl}^2\dot u^ju^k\partial_{l}
u^j-\partial_{k}\dot
u^j\partial_lu^k\partial_{l} u^j\right)dx\nonumber\\
&=&\mu \sigma^m\int \left(|\nabla\dot u|^2-\partial_k\dot
u^j\partial_ku^l\partial_l u^j+\partial_l\dot
u^j\partial_ku^k\partial_{l} u^j-\partial_k\dot
u^j\partial_lu^k\partial_{l}u^j\right) dx\nonumber\\
&\geq&\frac{7\mu}{8}\sigma^m\int |\nabla\dot u|^2dx-C \sigma^m\int
|\nabla u|^4 dx,\label{3.22}
\end{eqnarray}
where the Cauchy-Schwarz inequality was also used in the last
inequality. Analogously,
\begin{eqnarray}
L_3 &\geq&(\mu+\lambda) \sigma^m\int (\divg \dot u)^2dxds-C\
\sigma^m\int |\nabla \dot u||\nabla u|^2
dx\nonumber\\
&\geq& \sigma^m\int \left((\mu+\lambda)(\divg \dot
u)^2-\frac{\mu}{8}|\nabla \dot u|^2\right)dx-C \sigma^m\int |\nabla
u|^4 dx.\label{3.23}
\end{eqnarray}
In view of (\ref{3.16}), we obtain after integrating by parts that
\begin{eqnarray}
R_1 &=&- \sigma^m\int \divg \dot u\left(\divg(P(\rho) u)+(\rho
P'(\rho)-P(\rho))\divg
u\right)dx\nonumber\\
&&- \sigma^m\int \left(u^k\partial_k\divg\dot u+\partial_k\dot
u^j\partial_j u^k\right)\left(P(\rho)-P(\rho_s)\right)dx\nonumber\\
&=&-\sigma^m\int \divg \dot u\left[(\rho
P'(\rho)-P(\rho))\divg u +\divg \left(P(\rho_s)u\right)\right]dx\nonumber\\
&&-\sigma^m\int \partial_k\dot
u^j\partial_j u^k\left(P(\rho)-P(\rho_s)\right)dx\nonumber\\
&\leq& \frac{\mu}{8}\sigma^m\|\nabla\dot
u\|_{L^2}^2+C\sigma^m\left(\|\nabla
u\|_{L^2}^2+\|u\|_{L^6}^2\|\nabla\rho_s\|_{L^3}^2\right)\nonumber\\
&\leq&\frac{\mu}{8}\sigma^m\|\nabla\dot u\|_{L^2}^2+C\sigma^m
\|\nabla u\|_{L^2}^2,\label{3.24}
\end{eqnarray}
where we have also used (\ref{1.8}) and  (\ref{2.1}). Similarly, by
(\ref{1.1}) we have
\begin{eqnarray}
R_2&=&\sigma^m\int \rho u\cdot\nabla\left(\dot u^j\partial_j
f\right) dx- \sigma^m\int \partial_k\dot u^j\left(u^k(\rho-\rho_s)
\partial_jf\right)dx \nonumber\\
&\leq&C \sigma^m\left(\|\nabla \dot u\|_{L^2}\|\nabla
f\|_{L^3}\|u\|_{L^6}+\|\rho\dot u\|_{L^3}\|
u\|_{L^6}\|\nabla^2f\|_{L^2}\right) \nonumber\\
&\leq& \frac{\mu}{8}\sigma^m\|\nabla \dot u\|_{L^2}^2+C
\sigma^m\left(\|\rho^{1/2}\dot u\|_{L^2}^2+\|\nabla
u\|_{L^2}^2\right),\label{3.25}
\end{eqnarray}
where we have also used Cauchy-Schwarz inequality and the following
estimate:
\begin{eqnarray*}
\|\rho \dot u\|_{L^3}\leq C\|\rho^{1/2} \dot u\|_{L^2}^{1/2}\|\dot
u\|_{L^6}^{1/2}\leq C\|\rho^{1/2} \dot u\|_{L^2}^{1/2}\|\nabla\dot
u\|_{L^2}^{1/2}.
\end{eqnarray*}

Putting the estimates of (\ref{3.21})--(\ref{3.25}) into
(\ref{3.20}) gives
\begin{eqnarray}
&&\left(\sigma^m\|\rho^{1/2}\dot
u\|_{L^2}^2\right)_t+\mu\sigma^m\|\nabla
\dot u\|_{L^2}^2\nonumber\\
&&\quad\leq C\sigma^m\left(\|\nabla u\|_{L^2}^2+\|\nabla
u\|_{L^4}^4\right)+
C\left(m\sigma^{m-1}\sigma'+\sigma^m\right)\|\rho^{1/2}\dot
u\|_{L^2}^2.\label{3.26}
\end{eqnarray}
Thus, choosing $m=3$ in (\ref{3.26}) and integrating it over $(0,T)$
lead to (\ref{3.11}) immediately.\hfill$\square$

Next, we prove some important connections among the modified
``effective viscous flux", the density, and the gradient and
material derivative of the velocity,  which are crucial for our
further analysis. To do this, we need to make a full use of the
mathematical structure of steady state density $\rho_s$ of
(\ref{1.6}) and the standard $L^p$-estimate of elliptic system,
basing on a key observation due to Huang-Li-Xin \cite{HLX2006} (see
(\ref{3.29}) below).

\begin{lem}\label{lem3.3}Let $f\in H^2$ satisfy (\ref{1.7}) and $\rho_s=\rho_s(x)$ be the steady state
density of (\ref{1.6}). Assume that $(\rho,u)$ with $\rho\in [0,
2\tilde\rho]$ is a smooth solution of (\ref{1.1})--(\ref{1.5}) and
$F=F(x,t)$ is the modified ``effective viscous flux" defined by
(\ref{1.26}). Then there exists a positive constant $C$, depending
only on
$\bar\rho,\underline\rho,\tilde\rho,\inf\limits_{x\in\R^3}f(x)$ and
$\|f\|_{H^2}$, such that
\begin{eqnarray}
&&\|\nabla F\|_{L^2}+\|\nabla(\rho_s^{-1}{\rm curl} u)\|_{L^2} \leq
C\left(\|\rho\dot u\|_{L^2}+\|\nabla
u\|_{L^3}+\|\rho-\rho_s\|_{L^6}^2\right),\label{3.27}\\[2mm]
&&\|\nabla u\|_{L^6}\leq C\left(\|\rho\dot u\|_{L^2} +\|\nabla
u\|_{L^2}+\|\rho-\rho_s\|_{L^6}^2+\|\rho-\rho_s\|_{L^6}\right).
\label{3.28}
\end{eqnarray}
\end{lem}
\pf The proofs are analogous to those in \cite{Ho1995-1,HLX2010},
however, since  $\rho_s-\n_\infty\not\in H^3,$ we modify  the proof of $\|\nabla F\|_{L^2}$  slightly. Indeed, as observed in \cite{HLX2006}, one can utilize
(\ref{1.6}) to get that
\begin{eqnarray}
&&\rho_s^{-1}\left(\nabla P(\rho)-\rho\nabla
f\right)\nonumber\\
&&\quad=\rho_s^{-1}\left[\nabla\left(P(\rho)-P(\rho_s)\right)-\rho_s^{-1}(\rho-\rho_s)\nabla
P(\rho_s)\right] \nonumber\\
&&\quad=\nabla\left[\rho_s^{-1}\left(P(\rho)-P(\rho_s)\right)\right]
+\rho_s^{-2}\left[P(\rho)-P(\rho_s)-(\rho-\rho_s)
P'(\rho_s)\right]\nabla\rho_s.\label{3.29}
\end{eqnarray}
Thus, by virtue of (\ref{1.2}), (\ref{1.26}) and (\ref{3.29}), we
deduce
\begin{eqnarray}
&&
\rho_s^{-1}\rho\dot u-\nabla F+\mu{\rm curl}\left(\rho_s^{-1}{\rm
curl}
u\right) \nonumber \\
&&\quad    =-\left[(\lambda+2\mu)(\divg u)
\nabla\rho_s^{-1}-\mu\nabla
\rho_s^{-1}\times({\rm curl} u)\right] \nonumber \\
&&\qquad+\left[P(\rho)-P(\rho_s)-P'(\rho_s)(\rho-\rho_s)\right]
\nabla\rho_s^{-1} \nonumber \\
&&\quad  \triangleq G_1+G_2,\label{3.30}
\end{eqnarray}
which  gives
\begin{equation}
\Delta F=\div (\rho_s^{-1}\rho\dot u-G_1-G_2)\label{3.31}
\end{equation}
and
\begin{eqnarray}
\mu\Delta\left(\rho_s^{-1}{\rm curl}
u\right)&=&\mu \na \div \left(\rho_s^{-1}{\rm curl} u\right)-\mu
{\rm curl}{\rm curl}\left(\rho_s^{-1}{\rm curl} u\right)
\nonumber\\
&=&\mu \na \left( {\rm curl} u\cdot\na \rho_s^{-1}\right)-{\rm
curl}(\rho_s^{-1}\rho\dot u-G_1-G_2).\label{3.32}
\end{eqnarray}
Noticing that \bnn\ba
 \|G_1\|_{L^2}+\|G_2\|_{L^2} \leq& C\left(\|\nabla
\rho_s\|_{L^6}\|\nabla
u\|_{L^3}+\|\nabla\rho_s\|_{L^6}\|(\rho-\rho_s)^2\|_{L^3}\right)\\
 \leq& C\left(\|\nabla u\|_{L^3}+\|\rho-\rho_s\|_{L^6}^2\right),
\ea\enn we apply  the standard $L^2$-estimate   to
\eqref{3.31} and \eqref{3.32} to obtain
\begin{eqnarray*}
&&\|\nabla F\|_{L^2}+\|\nabla(\rho_s^{-1}{\rm curl} u)\|_{L^2} \nonumber\\
&&\quad\leq C\left( \|\rho\dot
u\|_{L^2}+\|G_1\|_{L^2}+\|G_2\|_{L^2}
+\|\nabla\rho_s\|_{L^6}\|\nabla u\|_{L^3}\right) \no &&\quad \leq
C\left(\|\rho\dot u\|_{L^2}+\|\nabla
u\|_{L^3}+\|\rho-\rho_s\|_{L^6}^2\right),
\end{eqnarray*} which gives (\ref{3.27}).

On the other hand, using (\ref{1.26}), (\ref{2.1}), and
(\ref{3.27}), we find
\begin{eqnarray*}
\|\nabla u\|_{L^6}&\leq& C\left(\|\divg u\|_{L^6}+\|{\rm curl}
u\|_{L^6}\right)\nonumber\\
&\leq& C \left(\|F\|_{L^6}+\|\rho-\rho_s\|_{L^6}+\|\rs^{-1}{\rm
curl} u\|_{L^6}\right)\nonumber\\
&\leq& C \left(\|\nabla F\|_{L^2}+\|\rho-\rho_s\|_{L^6}+\|\nabla
(\rs^{-1}{\rm curl}
u)\|_{L^2}\right) \nonumber\\
&\leq& C\left(\|\rho\dot u\|_{L^2}+\|\rho-\rho_s\|_{L^6}^2
+\|\rho-\rho_s\|_{L^6} +\|\nabla u\|_{L^2}\right)+\frac12\|\nabla
u\|_{L^6},
\end{eqnarray*}
where we have used the fact that $\|\nabla u\|_{L^3}^2\leq C\|\nabla
u\|_{L^2}\|\nabla u\|_{L^6}$ and Cauchy-Schwarz inequality in the
last inequality. This proves (\ref{3.28}), and thus, the proof of
Lemma \ref{lem3.3} is completed.\hfill$\square$

With the help of (\ref{3.4}) and Lemma \ref{lem3.3}, we can prove
\begin{lem}\label{lem3.4} Let $(\rho,u)$ with $\rho\in[0,2\tilde\rho]$ be a smooth
solution of (\ref{1.1})--(\ref{1.5}) on $\R^3\times (0,T]$. Then
there exists a positive constant $\varepsilon_0>0$, depending on
$\tilde\rho$, such that
\begin{equation}
 \int_0^T\sigma^3\left(\|\na u\|_{L^4}^4+\|\n-\rs\|_{L^4}^4+\|F\|_{L^4}^4+\|{\rm curl}u\|_{L^4}^4\right)dt\le CC_0.\label{3.34}
\end{equation}
provided $\Phi_1(T)+\Phi_2(T)\leq 2C_0^{1/2}$ and
$C_0\leq\varepsilon_0$.
\end{lem}
\pf Direct calculations from (\ref{1.1}) and  (\ref{1.26}) show
\begin{eqnarray*}
&&(\rho-\rho_s)_t+\frac{\rs}{2\mu+\lambda}\left(P(\rho)-P(\rho_s)\right)\nonumber\\
&&\quad=-u\cdot \na (\rho-\rho_s)-(\rho-\rho_s)\div u-u\cdot\na
\rho_s-\frac{\rho_s^2 F}{2\mu+\lambda}.
\end{eqnarray*}
Multiplying this by $4(\rho-\rho_s)^3$ in $L^2$ and integrating by
parts, we find
\begin{eqnarray}
&&\frac{d}{dt}\int (\rho-\rho_s)^4dx+
\frac{4}{2\mu+\lambda}\int\rho_s\left(P(\rho)-P(\rho_s)\right)(\rho-\rho_s)^3dx\nonumber\\
&&\quad\le C\int \left((\rho-\rho_s)^4|\na u|+|\rho-\rho_s|^3|u||\na
\rho_s|+ |\rho-\rho_s|^3|F|\right) dx\nonumber\\
&&\quad\le C\|\rho-\rho_s\|_{L^4}^2\left(\|\na
u\|_{L^2}+\|u\|_{L^6}\|\na\rs\|_{L^3}\right)+C\|\rho-\rho_s\|_{L^4}^3\|F\|_{L^4}\nonumber\\
&&\quad\le \de \|\rho-\rho_s\|_{L^4}^4+C(\de)\left(\|\na
u\|_{L^2}^2+\|F\|_{L^4}^4\right),\quad\delta>0,\label{3.35}
\end{eqnarray}
where we have used (\ref{1.8}), (\ref{2.1}) and Cauchy-Schwarz
inequality. Due to (\ref{1.8}), one has
\begin{eqnarray*} \rho_s\left(P(\rho)-P(\rho_s)\right)(\rho-\rho_s)^3&=&(\rho-\rho_s)^4\rs
\int_0^1P'(\alpha\n+(1-\al)\rho_s)d\al
\\
&\ge&  C(\rho-\rho_s)^4
\end{eqnarray*}
for some positive constant $C$ depending only on $A,
\gamma,\underline\rho,\bar\rho$ and $\tilde\rho$. Hence, choosing
$\de>0$ suitably small, multiplying (\ref{3.35}) by $\si^3$ and
integrating it over $(0,T)$, we infer from (\ref{3.7}) and
(\ref{3.9}) that
\begin{equation}
\int_0^T\si^3\|\n-\rs\|_{L^4}^4dt\le
CC_0+C\int_0^T\si^3\|F\|_{L^4}^4dt.\label{3.36}
\end{equation}

By (\ref{1.26}), we have
$$
\|\nabla u\|_{L^4}\leq C\left(\|\divg u\|_{L^4}+\|{\rm curl}
u\|_{L^4}\right)\leq C\left(\|F\|_{L^4}+\|\rho-\rho_s\|_{L^4}+\|\rs^{-1}{\rm
curl} u\|_{L^4}\right),
$$
which, combining with (\ref{3.36}) and (\ref{3.27}), yields that
\bnn\ba
& \int_0^T\sigma^3\left(\|\na u\|_{L^4}^4+\|\n-\rs\|_{L^4}^4+\|F\|_{L^4}^4 \right)dt \\
&  \le CC_0+C\int_0^T\sigma^3\left(\|F\|_{L^4}^4+\|\rs^{-1}{\rm curl}u\|_{L^4}^4\right)dt \\
& \leq CC_0+C\int_0^T\sigma^3\left(\|F\|_{L^2}\|\nabla
F\|_{L^2}^{3}+
\|\rs^{-1}{\rm curl}u\|_{L^2}\|\nabla(\rs^{-1}{\rm curl}u)\|_{L^2}^{3}\right)dt \\ &  \leq CC_0+C\int_0^T\sigma^3\left( \|\nabla u\|_{L^2}+C_0^{1/2}\right)
\left(\|\rho\dot u\|_{L^2}^{3}+\|\nabla u\|_{L^3}^3+\|\rho-\rho_s\|_{L^6}^6\right)dt\\ &  \leq CC_0+C\int_0^T\left( \sigma^{1/2}\|\nabla u\|_{L^2}+C_0^{1/2}\right)
\left(\sigma^{3}\|\rho\dot u\|_{L^2}^2\right)^{1/2}\si\|\rho\dot u\|_{L^2}^2dt \\ &  \quad +C\int_0^T  \|\nabla u\|_{L^2}
\left( \left( \sigma^{1/2}\|\nabla u\|_{L^2}\right)\sigma^{5/2}\|\nabla u\|_{L^4}^2 +\sigma^3 C_0^{1/2}\|\rho-\rho_s\|_{L^4}^2\right)dt \\ &   \quad+CC_0^{1/2}\int_0^T\sigma^3 \left( \|\nabla u\|_{L^2}^2+ \|\nabla u\|_{L^4}^4+\|\rho-\rho_s\|_{L^4}^4\right)dt \\&\le CC_0+CC_0^{1/2}\int_0^T\sigma^3 \left(  \|\nabla u\|_{L^4}^4+\|\rho-\rho_s\|_{L^4}^4\right)dt   ,\label{3.37}
\ea\enn where in the last inequality we have used
  $\Phi_1(T)+\Phi_2(T)\leq 2 C_0^{1/2}. $
This directly gives (\ref{3.34}) provided $C_0\le
\varepsilon_0\triangleq \min\{1,(2C)^{-2}\}$. \hfill$\square$

\begin{lem}\label{lem3.5} Let $(\rho,u)$ with
$\rho\in[0,2\tilde\rho]$ be a smooth solution of
(\ref{1.1})--(\ref{1.5}) on $\R^3\times(0,T]$. Then there exist
positive constants $K$ and $\varepsilon_1$, depending on
$\tilde\rho$ and $M$, such that
\begin{equation}
\Phi_3(\sigma(T))+\int_0^{\sigma(T)}\|\rho^{1/2}\dot u\|_{L^2}^2
dt\leq 2K,\label{3.38}
\end{equation}
provided $\Phi_3(\sigma(T))\leq 3K$ and $C_0\leq \varepsilon_1$.
\end{lem}
\pf Choosing $m=0$ in (\ref{3.19}) and integrating it over
$[0,\si(T)]$, we deduce from (\ref{3.7}), (\ref{3.9}), (\ref{3.28})
and Cauchy-Schwarz inequality that
\begin{eqnarray*}
&&\Phi_3(\sigma(T))+2\int_0^{\sigma(T)}\|\rho^{1/2} \dot u\|_{L^2}^2
dt\\
&&\quad\leq C\left(C_0+M\right)+C\int_0^{\sigma(T)}\|\nabla
u\|_{L^3}^3dt\\
&&\quad\leq C\left(C_0+M\right)+C\int_0^{\sigma(T)}\|\nabla
u\|_{L^2}^{3/2}\|\nabla
u\|_{L^6}^{3/2}dt\\
&&\quad\leq C\left(C_0+M\right) +C\int_0^{\sigma(T)}\|\nabla
u\|_{L^2}^{6}dt+\int_0^{\sigma(T)}\|\rho^{1/2}\dot
u\|_{L^2}^{2}dt\\
&&\quad\leq K+CC_0\left[\Phi_3(\sigma(T))\right]^2+
\int_0^{\sigma(T)}\|\rho^{1/2}\dot u\|_{L^2}^2dt
\end{eqnarray*}
with a positive constant $K\triangleq C(C_0+M)$ depending on
$\tilde\rho,C_0$ and $M$. As a result, we immediately obtain
(\ref{3.38}) provided $\Phi_3(\sigma(T))\leq 3K$ and
$C_0\leq\varepsilon_1\triangleq (9CK)^{-1}$.\hfill$\square$

\begin{lem}\label{lem3.6} Assume that $(\rho,u)$ satisfying
(\ref{3.4}) with $K>0$ as in Lemma \ref{lem3.5} is a smooth solution
of (\ref{1.1})--(\ref{1.5}) on $\R^3\times(0,T]$. Then,
 there exists a positive constant
$\varepsilon_2>0$, depending on $\tilde\rho$ and $M$, such that
\begin{equation}
\Phi_1(T)+\Phi_2(T)\leq C_0^{1/2},\label{3.39}
\end{equation} provided
$C_0\leq \varepsilon_2.$
\end{lem}
\pf By (\ref{3.4}), (\ref{3.9}) and (\ref{3.28}), we find
\begin{eqnarray*}
\si \|\na u\|^{2/3}_{L^6}&\le& C\si  \left(\|\rho\dot u\|_{L^2}
+\|\nabla
u\|_{L^2}+\|\rho-\rho_s\|_{L^6}^2+\|\rho-\rho_s\|_{L^6}\right)^{2/3}\\
&\le& C\left(C_0^{1/6}+C_0^{2/9}+C_0^{1/9}\right)\leq CC_0^{1/9}.
\end{eqnarray*}
Thus, using (\ref{3.7}), (\ref{3.9}), (\ref{3.28}), (\ref{3.34}) and
(\ref{3.38}), one derives from (\ref{3.10}) and (\ref{3.11}) that
\begin{eqnarray*}
&&\Phi_1(T)+\Phi_2(T)\\
&&\quad\le CC_0+C\int_{\sigma(T)}^T\sigma \|\na u\|_{L^3}^3dt+C\int_0^{\sigma(T)}\sigma \|\na u\|_{L^3}^3dt\\
&&\quad\le CC_0+C\int^T_{\si(T)}\left(\|\na
u\|_{L^2}^2+\sigma^3\|\na u\|_{L^4}^4\right)dt+C\int_0^{\si(T)}\si
\|\na u\|_{L^2}^{3/2}\|\na
u\|_{L^6}^{3/2}dt\\
&&\quad\le CC_0 +CC_0^{1/9}\int_0^{\si(T)}\|\nabla u\|_{L^2}^{7/6}
\left(\|\rho\dot u\|_{L^2}
+\|\nabla u\|_{L^2}+\|\rho-\rho_s\|_{L^6}^2+\|\rho-\rho_s\|_{L^6}\right)^{5/6}dt\\
&&\quad\le  CC_0 +CC_0^{2/3}\left(\int_0^{\si(T)}\left(\|\rho\dot
u\|_{L^2}^2 +\|\nabla
u\|_{L^2}^2+\|\rho-\rho_s\|_{L^6}^4+\|\rho-\rho_s\|_{L^6}^2\right)dt\right)^{5/12}\\
&&\quad\le  C_1C_0^{2/3}\leq C_0^{1/2},
\end{eqnarray*}
provided $C_0\le \ep_2\triangleq \min\{\ve_0,\ve_1,C_1^{-6}\}$.
The proof of Lemma \ref{lem3.6} is thus complete.\hfill$\square$

To complete the proof of Proposition \ref{pro3.1}, it remains to
prove the uniform upper bound of the density. To do this, we first
need the following refined estimates.
\begin{lem}\label{lem3.7}  Assume that $(\rho,u)$ satisfying
(\ref{3.4}) with $K>0$ as in Lemma \ref{lem3.5} is a smooth solution
of (\ref{1.1})--(\ref{1.5}) on $\R^3\times(0,T]$. Then
\begin{eqnarray}
&&\sup_{0\leq t\leq T}\|\nabla u\|_{L^2}^2+\int_0^T\|\rho^{1/2}\dot
u\|_{L^2}^2dt\leq
C,\label{3.40}\\
&&\sup_{0\leq t\leq T}\left(\sigma\| \rho^{1/2}\dot
u\|_{L^2}^2\right)+\int_0^T\sigma\|\nabla \dot u\|_{L^2}^2dt\leq
C,\label{3.41}
\end{eqnarray}
provided $C_0\leq\varepsilon_2$.
\end{lem}
\pf  We only have to prove (\ref{3.41})  since the estimate
(\ref{3.40}) is an immediate result of (\ref{3.38}) and
(\ref{3.39}). Choosing $m=1$ in (\ref{3.26}) and integrating it over
$(0,T)$ yields that
\begin{eqnarray*}
&&\sup_{0\leq t\leq T}\left(\sigma\|\rho^{1/2}\dot
u\|_{L^2}^2\right)+\int_0^T\sigma\|\nabla\dot u\|_{L^2}^2dt\nonumber\\
&&\quad\leq C\left(C_0+\int_0^T\|\rho^{1/2}\dot
u\|_{L^2}^2dt+\int_0^T\sigma\|\nabla u\|_{L^4}^4dt\right)\nonumber\\
&&\quad\leq C+C\int_0^{\sigma(T)}\sigma\|\nabla u\|_{L^2}\|\nabla
u\|_{L^6}^{3}dt
\nonumber\\
&&\quad\leq C+C\int_0^{\sigma(T)}\sigma \left(\|\rho\dot u\|_{L^2}^3
+\|\nabla u\|_{L^2}^3+\|\rho-\rho_s\|_{L^6}^6 +\|\rho-\rho_s\|_{L^6}^3 \right)dt\\
&&\quad\leq C+C\sup_{0\leq t\leq
T}\left(\sigma^{1/2}\|\rho^{1/2}\dot u\|_{L^2}\right),
\end{eqnarray*}
where we have used (\ref{3.9}), (\ref{3.34}), (\ref{3.28}) and
(\ref{3.40}). Combining this with Young inequality immediately leads
to (\ref{3.41}).\hfill$\square$

We are now ready to prove the uniform upper bound of the density,
which is in fact the key to extend the local smooth solution to be
global and   will be proved by modifying the arguments in
\cite{HLX2006,LX2006,HLX2010}  basing on the Zlotnik inequality (cf.
Lemma \ref{lem2.2}) and the standard ``effective viscous flux"
$\tilde F$ (see (\ref{1.27})) in a similar form as the one in
\cite{Li1998,Ho1995-1,HLX2010}.

\begin{lem}\label{lem3.8}  Assume that $(\rho,u)$ satisfying
(\ref{3.4}) with $K>0$ as in Lemma \ref{lem3.5} is a smooth solution
of (\ref{1.1})--(\ref{1.5}) on $\R^3\times(0,T]$. Then there exists
a positive constant $\tilde\varepsilon>0$, depending on $\tilde\rho$
and $M$, such that for all $  (x,t)\in\R^3\times[0,T],$
$$
\rho(x,t)\leq \frac{7}{4}\tilde\rho
$$
provided $C_0\leq \tilde\varepsilon$.
\end{lem}
\pf Let  $\tilde F $ be the standard
``effective viscous flux" defined by  (\ref{1.27}). Then (\ref{1.1})
can be written in the form:
$$
D_t\rho=g(\rho)+b'(t),
$$
where
$$D_t\triangleq
\partial_t+u\cdot\nabla,\quad
g(\rho)\triangleq-\frac{A\rho}{\lambda+2\mu}\left(\rho^\gamma-\rho_s^\gamma\right),\quad
b(t)\triangleq-\frac{1}{\lambda+2\mu}\int_0^t\rho \tilde F d\tau.
$$

In order to apply Lemma \ref{lem2.2}, we have to deal with $b(t)$ since $\lim\limits_{\n\rightarrow +\infty}g(\n)=-\infty.$
To do this, we first observe from (\ref{1.2}) and (\ref{1.6}) that
\begin{equation}
\la{3.42} \Delta \tilde F=\divg(\rho\dot u-(\rho-\rho_s)\nabla f).
\end{equation}
Applying the standard $L^p$-estimate to the elliptic problem
(\ref{3.42}) gives \be\ba \|\na \tilde F\|_{L^4}&\le  C\|\n\dot
u\|_{L^4}+C\|(\rho-\rho_s)\nabla f\|_{L^4} \\
&\le  C\|\rho\dot u\|_{L^2}^{1/4}\|\rho\dot
u\|_{L^6}^{3/4}+C\|\rho-\rho_s\|_{L^{12}}\|\na f\|_{L^6} \\
&\le  C\si^{-1/8} \|\na \dot
u\|_{L^2}^{3/4}+CC_0^{1/12},\label{3.43} \ea\ee where we have also
used (\ref{1.8}), (\ref{3.9}) and (\ref{3.41}). It thus follows from
(\ref{1.27}), (\ref{2.2}), (\ref{3.9}), (\ref{3.40}) and
(\ref{3.43}) that
\begin{eqnarray}
\|\tilde F\|_{L^\infty}&\leq& C\|\tilde F\|_{L^2}^{1/7}\|\nabla \tilde F\|_{L^4}^{6/7} \nonumber\\
&\leq& C\left(\|\nabla
u\|_{L^2}^{1/7}+\|\rho-\rho_s\|_{L^2}^{1/7}\right) \left(\si^{-1/8}
\|\na \dot u\|_{L^2}^{3/4}+ C_0^{1/12}
\right)^{6/7} \nonumber\\
&\leq&C  \si^{-3/28} \|\na \dot u\|_{L^2}^{9/14}+CC_0^{1/14}.
\label{3.44}
\end{eqnarray}

For $0\leq t_1<t_2\leq\sigma(T)\leq1$, we deduce from (\ref{3.39}),
(\ref{3.41}) and (\ref{3.44}) that if  $C_0\leq \varepsilon_2$, then
\begin{eqnarray*}
&&|b(t_2)-b(t_1)|\nonumber\\
&&\quad\leq C\int_0^{\sigma(T)}\|\tilde F\|_{L^\infty}dt\nonumber\\
&&\quad\leq C C_0^{1/14}+C
\int_0^{\sigma(T)}\sigma^{-4/7}\left(\sigma\|\nabla \dot
u\|_{L^2}^2\right)^{1/4}\left(\sigma^3\|\nabla \dot
u\|_{L^2}^2\right)^{1/14}dt \\
&&\quad\leq C C_0^{1/14}+C \left(\int_0^{\sigma(T)}\sigma\|\nabla
\dot
u\|_{L^2}^2dt\right)^{1/4}\left(\int_0^{\sigma(T)}\sigma^3\|\nabla
\dot u\|_{L^2}^2dt\right)^{1/14}\\
&&\quad\leq C C_0^{1/14}+C\Phi_2^{1/14}(T)\leq C C_0^{1/28}.
\end{eqnarray*}
Therefore, for any $t\in[0,\sigma(T)]$, one can choose $N_0, N_1$ in
(\ref{2.3}) and $\xi^*$ in (\ref{2.5}) as follows:
$$
N_0= C C_0^{1/28},\quad N_1= 0,\quad \xi^*=\bar\rho.
$$
Then, due to the facts that $\underline\rho\leq\rho_s\leq\bar\rho$
(see Proposition \ref{pro1.1}) and
$$
g(\xi)\leq
-\frac{A\xi}{\lambda+2\mu}\left(\xi^\gamma-\bar\rho^\gamma\right)\leq
-N_1= 0,\quad\forall\;\xi\geq\xi^*=\bar\rho,
$$
it follows from (\ref{2.4}) that (keeping in mind that
$0\leq\rho_0\leq\tilde\rho$ and $\tilde\rho\geq \bar\rho+1$)
\begin{equation}
\sup_{0\leq t\leq \sigma(T)}\|\rho(t)\|_{L^\infty}\leq
\max\{\tilde\rho,\bar\rho\}+N_0\leq\tilde\rho+CC_0^{1/28}\leq\frac{3}{2}\tilde\rho,\label{3.45}
\end{equation}
provided the initial energy $C_0$ is chosen to be such that
$$
C_0\leq\min\{\varepsilon_2,\varepsilon_3\}\quad{\rm
with}\quad\varepsilon_3\triangleq\left(\frac{\tilde\rho}{2C}\right)^{28}.
$$

For any $\sigma(T)\leq t_1<t_2\leq T$, we have from (\ref{3.9}),
(\ref{3.39}) and (\ref{3.40}) that
\begin{eqnarray*}
&&|b(t_2)-b(t_1)|\\
&&\quad\leq C\int_{t_1}^{t_2}\|\tilde F\|_{L^\infty}dt \leq
C\int_{t_1}^{t_2}\|\tilde F\|_{L^2}^{1/7}\|\nabla\tilde
F\|_{L^4}^{6/7}dt\nonumber\\
&&\quad\leq C\int_{t_1}^{t_2}\left(\|\rho\dot
u\|_{L^2}^{1/4}\|\rho\dot
u\|_{L^6}^{3/4}+\|\rho-\rho_s\|_{L^{12}}\|\na f\|_{L^6}\right)^{6/7}dt\\
&&\quad\leq CC_0^{1/14}(t_2-t_1) +C\int_{t_1}^{t_2}\|\n\dot
u\|_{L^2}^{3/14}\|\dot
u\|_{L^6}^{9/14}dt\nonumber\\
&&\quad \leq\left( C
C_0^{1/14}+\frac{A}{2(2\mu+\lambda)}\right)(t_2-t_1)
+C\int_{\sigma(T)}^{T}\sigma^3\left(\|\rho^{1/2}\dot u\|_{L^2}^2+
\|\rho \dot u\|_{L^2}^2\right)dt \nonumber\\
&&\quad \leq \frac{A}{2\mu+\lambda} (t_2-t_1)+CC_0^{1/2},
\end{eqnarray*}
where in the last inequality we have chosen $C_0$ to be such that
$$
C_0\le \min\{\ve_2,\ve_3,\ve_4\}\quad{\rm with}\quad \ve_4\triangleq
\left(\frac{A}{2C(2\mu+\lambda)}\right)^{14}.
$$
Therefore, for any $t\in[\sigma(T),T]$, we can choose $N_0$ and $N_1
$ in (\ref{2.3}) and $\xi^*$ in (\ref{2.5}) as follows:
$$
N_0=C C_0^{1/2},\quad
N_1=\frac{A}{2\mu+\lambda},\quad\xi^*=\bar\rho+1.
$$
Then it is easy to check that for any $\xi\geq\xi^*,$
$$
g(\xi)\leq -\frac{A\xi}{2\mu+\lambda}(\xi^\gamma-\bar\rho^\gamma)
\leq-N_1=-\frac{A}{2\mu+\lambda},
$$
and hence, one has from (\ref{2.4}) and (\ref{3.45}) that
\begin{equation}
\sup_{\sigma(T)\leq t\leq
T}\|\rho(t)\|_{L^\infty}\leq\max\left\{\frac{3}{2}\tilde\rho,\bar\rho+1\right\}+N_0\leq
\frac{3}{2}\tilde\rho+C
C_0^{1/2}\leq\frac{7}{4}\tilde\rho,\label{3.46}
\end{equation}
provided the initial energy $C_0$ satisfies
$$
C_0\leq\tilde\varepsilon\triangleq\min\{\varepsilon_2,\varepsilon_3,
\varepsilon_4,\ve_5\}\quad{\rm
with}\quad\varepsilon_5\triangleq\left(\frac{\tilde\rho}{4C}\right)^2.
$$
The combination of (\ref{3.45}) with (\ref{3.46}) proves Lemma
\ref{lem3.8}.\hfill$\square$

\section{Time-dependent higher-order estimates}
In this section, we prove the higher-order estimates of the smooth
solution $(\rho,u)$ to \eqref{1.1}-\eqref{1.5}, which are needed for
the existence of classical solutions. From now on, we always assume
that the initial energy $C_0$ satisfies \eqref{3.6}  and that the
conditions of Theorem \ref{thm1.1} are satisfied. We also denote by
$C$ the various positive constants which may depend on
$\rho_0,u_0,\|g\|_{L^2},f,$ $\lambda,$ $\mu,$ $A,$ $\gamma,$
$\underline\rho,$ $\bar\rho,\tilde\rho,M,$ and $T$ as well.

We begin with the $L^2$-estimate on the material derivative of the velocity.
\begin{lem}\label{lem4.1}For any given $T>0$, there exists a positive constant
$C(T)$ such that
\begin{equation}
\sup_{0\leq t\leq T}\|\rho^{1/2}\dot u\|_{L^2}^2 +\int_0^T\|\nabla
\dot u\|_{L^2}^2 dt\leq C(T).\label{4.1}
\end{equation}
\end{lem}
\pf Choosing $m=0$ in (\ref{3.26}) and integrating it over $(0,T)$,
we have
\begin{eqnarray*}
&&\sup_{0\leq t\leq T}\int \rho|\dot u|^2dx- \int \rho|\dot
u|^2(x,0)dx+\int_0^T\|\nabla\dot
u\|_{L^2}^2dt\\
&&\quad\leq C\int_0^T\left(\|\nabla u\|_{L^2}^2+\|\rho^{1/2}\dot
u\|_{L^2}^2\right)dt+C\int_0^T\|\nabla u\|_{L^4}^4dt\nonumber\\
&&\quad\leq C+C\int_0^T\|\nabla u\|_{L^2}\|\nabla u\|_{L^6}^3dt\nonumber\\
&&\quad\leq C+C\int_0^T\|\rho\dot u\|_{L^2}^3 dt\nonumber\\
&&\quad\leq C+C\sup_{0\leq t\leq T}\|\rho^{1/2}\dot u\|_{L^2},
\end{eqnarray*}
where we have used (\ref{3.7}), (\ref{3.9}), (\ref{3.28}) and
(\ref{3.40}). Combining this with Young inequality gives
(\ref{4.1}), since the compatibility condition (\ref{1.14}) implies
that $\sqrt\rho\dot u|_{t=0}=\sqrt{\rho_0}(\nabla f-g)\in L^2$ is
well defined.\hfill$\square$

Next, similar to that in \cite{HLX2010,HLX2010-2}, we utilize the
Beale-Kato-Majda-type inequality (see Lemma \ref{lem2.3}) to prove
the important estimates on the gradients of $(\rho,u)$.
\begin{lem}\label{lem4.2}There exists a positive constant
$C=C(T)$ such that
\begin{equation}
\sup_{0\leq t\leq T}\left(\|\nabla\rho\|_{L^2{\cap}L^6}+\|\nabla
u\|_{H^1}\right)+\int_0^T\|\nabla u\|_{L^\infty}dt\leq
C(T).\label{4.2}
\end{equation}
\end{lem}
\pf  Since $\mathcal{L}=-\mu\Delta-(\mu+\lambda)\nabla\divg$ is a
strong elliptic operator (see \cite{CCK2004} for instance), applying
the $L^p$-estimate of elliptic system to (\ref{1.2}) gives that for
$ 2\le p\le 6,$
\begin{equation}
\|\nabla^2 u\|_{L^p}\leq C\left(\|\rho\dot
u\|_{L^p}+\|\nabla\rho\|_{L^p}+\|\nabla f\|_{L^p}\right).
\label{4.3}
\end{equation}
By integration by parts, we easily derive from (\ref{1.1}) that
\begin{eqnarray}
\frac{d}{dt}\|\nabla\rho\|_{L^p}&\leq& C\|\nabla
u\|_{L^\infty}\|\nabla\rho\|_{L^p}+C\|\nabla^2u\|_{L^p}
\nonumber \\
&\leq& C\left(1+\|\nabla
u\|_{L^\infty}\right)\|\nabla\rho\|_{L^p}+C\left(1+\|\rho\dot
u\|_{L^p}\right),\label{4.4}
\end{eqnarray}
where (\ref{4.3}) was used in the second inequality. Using
(\ref{3.40}), (\ref{4.3}) and (\ref{2.1}), we deduce from
(\ref{2.6}) that
\begin{eqnarray}
\|\nabla u\|_{L^\infty} &\leq& C\left(\|{\rm div}
u\|_{L^\infty}+\|{\rm curl}
u\|_{L^\infty}\right)\ln\left(e+\|\rho\dot
u\|_{L^6}+\|\nabla\rho\|_{L^6}\right)+C\nonumber\\
&\leq& C+C\left(\|{\rm div} u\|_{L^\infty}+\|{\rm curl}
u\|_{L^\infty}\right)\ln\left(e+\|\nabla\dot
u\|_{L^2}\right)\nonumber\\
&&+ C\left(\|{\rm div} u\|_{L^\infty}+\|{\rm curl}
u\|_{L^\infty}\right)\ln\left(e+\|\nabla\rho\|_{L^6}\right).\label{4.5}
\end{eqnarray}

Set
$$
\Phi(t)\triangleq e+\|\nabla\rho\|_{L^6},\quad
\Psi(t)\triangleq1+\left(\|{\rm div} u\|_{L^\infty}+\|{\rm curl}
u\|_{L^\infty}\right)\ln\left(e+\|\nabla\dot
u\|_{L^2}\right)+\|\nabla\dot u\|_{L^2}.
$$
Then it follows from (\ref{4.4}) with $p=6$ and (\ref{4.5}) that
$$
\Phi'(t)\leq C\Psi(t)\Phi(t)\ln\Phi(t),
$$
which, together with the fact that $\Phi(t)>1$, implies
\begin{equation}
\frac{d}{dt}\ln\Phi(t)\leq  C\Psi(t)\ln\Phi(t).\label{4.6}
\end{equation}
Applying the standard $L^p$-estimate to (\ref{3.31}) and
(\ref{3.32}) yields that
\begin{eqnarray*}
\|\nabla F\|_{L^6}+\|\nabla (\rs^{-1} {\rm
curl}u)\|_{L^6}&\leq&C\left(\|\rho\dot
u\|_{L^6}+\|G_1\|_{L^6}+\|G_2\|_{L^6}+\|\nabla
u\nabla\rho_s^{-1}\|_{L^6}\right)\nonumber\\
&\leq&C\left(\|\nabla\dot
u\|_{L^2}+\|\nabla\rho_s\|_{L^\infty}\|\nabla
u\|_{L^6}+\|\nabla\rho_s\|_{L^6}\right)\nonumber\\
&\leq&C\left(\|\nabla\dot u\|_{L^2}+1\right),
\end{eqnarray*}
where we have also used (\ref{1.9}), (\ref{3.9}), (\ref{3.28}),
(\ref{3.40}), (\ref{4.1}) and Lemma \ref{lem2.1}. Thus,
\begin{eqnarray}
\int_0^T\Psi(t)dt&\leq& C+C\int_0^T\left(\|\nabla \dot
u\|_{L^2}^2+\|\divg u\|_{L^\infty}^2+\|{\rm curl}
u\|_{L^\infty}^2\right)dt\nonumber\\
&\leq&C+C\int_0^T\left(\|P(\rho)-P(\rho_s)\|_{L^\infty}^2+\|F\|_{L^\infty}^2
+\|\rs^{-1}{\rm
curl}u\|_{L^\infty}^2\right)dt\nonumber\\
&\leq&C+C\int_0^T\left(\|F\|_{L^2}^2+\|{\rm
curl}u\|_{L^2}^2+\|\nabla F\|_{L^6}^{2}+\|\nabla(\rs^{-1}{\rm
curl}u)\|_{L^6}^{2}\right)dt\nonumber\\
&\leq&C+C\int_0^T\|\nabla\dot u\|_{L^2}^2dt\leq C,\label{4.7}
\end{eqnarray}
and consequently, it follows from (\ref{4.6}) and Gronwall's
inequality that
\begin{equation}
\sup_{0\leq t\leq T}\|\nabla\rho(t)\|_{L^6}\leq\sup_{0\leq t\leq
T}\Phi(t)\leq C.\label{4.8}
\end{equation}
As a result of (\ref{4.5}), (\ref{4.7}) and (\ref{4.8}), we obtain
\begin{equation}
\int_0^T\|\nabla u\|_{L^\infty}dt\leq  C+C\int_0^T\Psi(t)dt\leq
C.\label{4.9}
\end{equation}
Using (\ref{4.1}) and (\ref{4.9}), we infer from (\ref{4.4}) with
$p=2$ and Gronwall's inequality that
$$
\sup_{0\leq t\leq T}\|\nabla\rho(t)\|_{L^2}\leq C,
$$
which, together with (\ref{4.1}) and (\ref{4.3}), also implies that
$\|\nabla u\|_{H^1}\leq C$. The proof of Lemma \ref{lem4.2} is
therefore completed.\hfill$\square$

By virtue of Lemmas \ref{lem4.1} and \ref{lem4.2}, we easily obtain
\begin{lem}\label{lem4.3}For any given $T>0$, it holds that
\begin{eqnarray}
&&\sup_{0\leq t\leq T}\|\rho^{1/2}u_t\|_{L^2}^2+\int_0^T\|\nabla
u_t\|_{L^2}^2dt\leq C(T),\label{4.10}\\
&&\sup_{0\leq t\leq T}\left(\|\nabla \rho\|_{H^1}+\|\nabla
P\|_{H^1}\right)+\int_0^T\|\nabla^2 u\|_{H^1}^2dt\leq
C(T),\label{4.11}\\
&&\sup_{0\leq t\leq
T}\left(\|\rho_t\|_{H^1}+\|P_t\|_{H^1}\right)+\int_0^T\left(\|\rho_{tt}\|_{L^2}^2+\|P_{tt}\|_{L^2}^2\right)dt\leq
C(T).\label{4.12}
\end{eqnarray}
\end{lem}
\pf First, (\ref{4.10}) follows directly from Lemmas \ref{lem4.1}
and \ref{lem4.2}.

Next, we prove (\ref{4.11}). Since $P(\rho)=A\rho^\gamma$ satisfies
\begin{equation}
P_t+u\cdot\nabla P+\gamma P\divg u=0,  \label{4.13}
\end{equation}
from which and (\ref{1.1}) we have by (\ref{4.2}) and direct
computations that
\begin{eqnarray}
&&\frac{d}{dt}\left(\|\nabla^2\rho\|_{L^2}^2+\|\nabla^2
P\|_{L^2}^2\right)\nonumber\\
&&\quad\leq C\|\nabla
u\|_{L^\infty}\left(\|\nabla^2\rho\|_{L^2}^2+\|\nabla^2
P\|_{L^2}^2\right)+C\|\nabla^3
u\|_{L^2}\left(\|\nabla^2\rho\|_{L^2}+\|\nabla^2
P\|_{L^2}\right)\nonumber\\
&&\qquad+C\|\nabla^2 u\|_{L^3}\left(\|\nabla\rho\|_{L^6}+\|\nabla
P\|_{L^6}\right)\left(\|\nabla^2\rho\|_{L^2}+\|\nabla^2
P\|_{L^2}\right)\nonumber\\
&&\quad\leq C\left(1+\|\nabla
u\|_{L^\infty}\right)\left(\|\nabla^2\rho\|_{L^2}^2+\|\nabla
P\|_{L^2}^2\right)+C\|\nabla^2 u\|_{H^1}^2.\label{4.14}
\end{eqnarray}
Using (\ref{4.2}), (\ref{4.10}) and Lemma \ref{lem2.1}, we deduce
from (\ref{1.2}) and the standard $L^2$-estimate of elliptic system
that
\begin{eqnarray}
\|\nabla^2 u\|_{H^1}&\leq& C\left(\|\rho u_t\|_{H^1}+\|\rho
u\cdot\nabla u\|_{H^1}+\|\nabla
P\|_{H^1}+\|\rho\nabla f\|_{H^1}\right)\nonumber\\
&\leq& C\left(1+\|\nabla
u_t\|_{L^2}+\|\nabla\rho\|_{L^3}\|u_t\|_{L^6}+\|\nabla^2
P\|_{L^2}\right)\nonumber\\
&\leq&C\left(1+\|\nabla
u_t\|_{L^2}+\|\nabla^2P\|_{L^2}\right).\label{4.15}
\end{eqnarray}
Thus, putting (\ref{4.15}) into (\ref{4.14}) and using Gronwall
inequality, we immediately arrive at (\ref{4.11}) since it holds
that $\|\nabla u\|_{L^\infty}+\|\nabla u_t\|_{L^2}^2\in L^1(0,T)$
due to (\ref{4.2}) and (\ref{4.10}) .

Finally, we prove (\ref{4.12}).  Thanks to (\ref{4.2}) and
(\ref{4.11}), it follows from (\ref{1.1}) and (\ref{4.13}) that
\begin{eqnarray}
\|\rho_t\|_{L^2}+\|P_t\|_{L^2}&\leq&
C\|u\|_{L^\infty}\left(\|\nabla\rho\|_{L^2}+\|\nabla
P\|_{L^2}\right)+C\|\nabla u\|_{L^2}\leq C,\label{4.16}
\\[2mm]
\|\nabla\rho_t\|_{L^2}+\|\nabla P_t\|_{L^2}&\leq& C\|\nabla^2
u\|_{L^2}+C\|u\|_{L^\infty}\left(\|\nabla^2\rho\|_{L^2}+\|\nabla^2P\|_{L^2}\right)\nonumber\\
&&+C\|\nabla u\|_{L^3}\left(\|\nabla\rho\|_{L^6}+\|\nabla
P\|_{L^6}\right)\leq C,\label{4.17}
\end{eqnarray}
where Sobolev inequalities were used to get that
$\|u\|_{L^\infty}+\|\nabla u\|_{L^3}\leq C\|\nabla u\|_{H^1}\leq C$.
Moreover, since (\ref{4.13}) implies
$$
P_{tt}+u_t\cdot\nabla P+u\cdot\nabla P_t+\gamma P_t\divg u+\gamma
P\divg u_t=0,
$$
one obtains after using (\ref{4.2}), (\ref{4.10}), (\ref{4.11}),
(\ref{4.16}), and (\ref{4.17}) that
\begin{eqnarray}
\int_0^T\|P_{tt}\|_{L^2}^2dt&\leq&
C\int_0^T\left(\|u_t\|_{L^6}\|\nabla P\|_{L^3}+\|\nabla
P_t\|_{L^2}+\|\nabla u\|_{L^\infty}\|P_t\|_{L^2}+\|\nabla
u_t\|_{L^2}\right)^2dt\nonumber\\
&\leq& C+C\int_0^T\left(\|\nabla u\|_{H^2}^2+\|\nabla
u_t\|_{L^2}^2\right)dt\leq C.\label{4.18}
\end{eqnarray}
In the same way, one also has $\|\rho_{tt}\|_{L^2}\in L^2(0,T)$. So,
combining this with (\ref{4.16})--(\ref{4.18}) completes the proof
of (\ref{4.12}).\hfill$\square$

In order to prove the solution obtained is indeed a classical one on
the time-interval $[\tau, T]$ for any $0<\tau<T<\infty$, we need
some further estimates on the higher order derivatives of $(\rho,
u)$. However, due to the weaker compatibility condition (\ref{1.14})
(cf. (\ref{1.19}), (\ref{1.20})), the methods used in \cite{HLX2010}
cannot be applied any more. To overcome this difficulty, we need the
following initial-layer analysis.
\begin{lem}\label{lem4.4}Let $\sigma\triangleq\min\{1,t\}$. Then it holds for any given $T>0$ that
\begin{eqnarray}
 \sup_{0\leq t\leq T} \sigma\left(\|\nabla
u\|_{H^2}^2+\|\nabla u_t\|_{L^2}^2\right) +\int_0^T\sigma
\left(\|\rho^{1/2}u_{tt}\|_{L^2}^2+\|\nabla
u_t\|_{H^1}^2\right)dt\leq C(T).\label{4.19}
\end{eqnarray}
\end{lem}
\pf First, differentiating (\ref{1.2}) with respect to $t$ gives
\be \la{a4.20}\ba&  \mu\Delta u_t+(\mu+\lambda)\nabla{\rm div}u_t\\ &=\n u_{tt}-\rho_t \nabla f +\rho_t u_t+\rho_tu\cdot\nabla
u+\rho u_t\cdot\nabla u+\rho u\cdot\nabla
u_t+\nabla P_t.\ea\ee
Multiplying (\ref{a4.20})  by $u_{tt}$ and integrating the resulting equality over $\r^3$ yield
\begin{eqnarray}
&&\frac{1}{2}\frac{d}{dt}\int \left[\mu|\nabla
u_t|^2+(\mu+\lambda)(\divg
u_t)^2\right]dx+\int \rho|u_{tt}|^2dx\nonumber\\
&&\quad=\int \left(\rho_t \nabla f -\rho_t u_t-\rho_tu\cdot\nabla
u-\rho u_t\cdot\nabla u-\rho u\cdot\nabla
u_t-\nabla P_t\right)\cdot u_{tt}dx\nonumber\\
&&\quad=\frac{d}{dt}\int \left[-\frac{1}{2}\rho_t
|u_t|^2+\left(\rho_t \nabla f
-\rho_tu\cdot\nabla u-\nabla P_t\right)\cdot u_{t}\right]dx\nonumber\\
&&\qquad-\int \left(\rho_{tt} \nabla f-\rho_{tt}u\cdot\nabla u-\rho_tu_t\cdot\nabla u-\rho_tu\cdot\nabla u_t\right)\cdot u_{t}dx\nonumber\\
&&\qquad+\frac{1}{2}\int \rho_{tt}|u_t|^2dx-\int P_{tt}\divg u_t
dx-\int \left(\rho u_t\cdot\nabla u+\rho u\cdot\nabla
u_t\right)\cdot u_{tt}dx\nonumber\\
&&\quad\triangleq\frac{d}{dt}I_0+\sum_{i=1}^4I_i.\label{4.20}
\end{eqnarray}

Then, we estimate each term on the right-hand side of (\ref{4.20}).
Using (\ref{1.1}) and integrating by parts, we see from
(\ref{4.2}) and (\ref{4.10})--(\ref{4.12}) that
\begin{eqnarray*}
I_0&=&\int \left[-\rho u\cdot\nabla u_t \cdot u_t
+\rho_t\left(\nabla f -u\cdot\nabla
u\right)\cdot u_t+P_t\divg u_t\right] dx\nonumber\\
&\leq&C\|u\|_{L^\infty}\|\rho^{1/2} u_t\|_{L^2}\|\nabla
u_t\|_{L^2}+C\|\rho_t\|_{L^2}\|\nabla
f\|_{L^3}\|u_t\|_{L^6}\nonumber\\
&&+C\|\rho_t\|_{L^2}\|u\|_{L^\infty}\|\nabla
u\|_{L^3}\|u_t\|_{L^6}+C\|P_t\|_{L^2}\|\nabla
u_t\|_{L^2}\nonumber\\
&\leq&\frac{\mu}{4}\|\nabla u_t\|_{L^2}^2+C.
\end{eqnarray*}
Using (\ref{4.2}) and (\ref{4.10})--(\ref{4.12}) again, one infers
that
\begin{eqnarray*}
I_1&\leq& C\|\rho_{tt}\|_{L^2}\|\nabla
f\|_{L^3}\|u_t\|_{L^6}+C\|\rho_{tt}\|_{L^2}\|u\|_{L^\infty}\|\nabla
u\|_{L^3}\|u_t\|_{L^6}\nonumber\\
&&+C\|\rho_t\|_{L^2}\|u_t\|_{L^6}^2\|\nabla
u\|_{L^6}+C\|\rho_t\|_{L^6}\|u\|_{L^6}\|\nabla
u_t\|_{L^2}\|u_t\|_{L^6}\nonumber\\
&\leq& C\|\rho_{tt}\|_{L^2}\|\nabla u_t\|_{L^2}+C\|\nabla
u_t\|_{L^2}^2\leq C\left(\|\rho_{tt}\|_{L^2}^2+\|\nabla
u_t\|_{L^2}^2\right).
\end{eqnarray*}
By (\ref{1.1}) and integration by parts, we have from (\ref{4.10})
and (\ref{4.12}) that
\begin{eqnarray*}
I_2&=&\int (\rho u)_t\cdot\nabla u_t\cdot u_tdx=\int \left(\rho_t
u\cdot\nabla u_t\cdot u_t+\rho
u_t\cdot\nabla u_t\cdot u_t \right)dx\nonumber\\
&\leq&C\|\rho_t\|_{L^3}\|u\|_{L^\infty}\|\nabla
u_t\|_{L^2}\|u_t\|_{L^6}+C\|\rho u_t\|_{L^3}\|\nabla
u_t\|_{L^2}\|u_t\|_{L^6}\nonumber\\
&\leq&C\|\nabla u_t\|_{L^2}^2+C\|\rho u_t\|_{L^2}^{1/2}\|
u_t\|_{L^6}^{1/2}\|\nabla u_t\|_{L^2}^2\leq C\left(1+\|\nabla
u_t\|_{L^2}^4\right).
\end{eqnarray*}
Obviously, $I_3\leq C\left(\|P_{tt}\|_{L^2}^2+\|\nabla
u_t\|_{L^2}^2\right)$. Finally, it follows from (\ref{4.2}) and
Lemma \ref{lem2.1} that
\begin{eqnarray*}
I_4&\leq&C\left(\|u_t\|_{L^6}\|\nabla
u\|_{L^3}+\|u\|_{L^\infty}\|\nabla u_t\|_{L^2}\right)\|\rho^{1/2}
u_{tt}\|_{L^2}\\
&\leq& \frac{1}{2}\|\rho^{1/2} u_{tt}\|_{L^2}^2+C\|\nabla
u_t\|_{L^2}^2.
\end{eqnarray*}

Thus, putting the estimates of $I_i$ into (\ref{4.20}), multiplying
the resulting inequality by $\sigma(t)$ and integrating it over
$[0,T]$, we infer from (\ref{4.10}), (\ref{4.12}) and Gronwall's
inequality that
$$
\sup_{0\leq t\leq T}\left(\sigma\|\nabla
u_t\|_{L^2}^2\right)+\int_{0}^T\sigma\|\rho^{1/2} u_{tt}\|_{L^2}^2
dt\leq C,
$$
which, together with (\ref{4.11}) and (\ref{4.15}), gives
\begin{equation}
\sup_{0\leq t\leq T} \sigma\left(\|\nabla u\|_{H^2}^2+\|\nabla
u_t\|_{L^2}^2\right) +\int_0^T\sigma\|
\rho^{1/2}u_{tt}\|_{L^2}^2dt\leq C.\label{4.21}
\end{equation}

 Finally, applying the standard $L^2$-estimate to the elliptic
system \eqref{a4.20} together with  (\ref{4.2}) and (\ref{4.12}) gives
\begin{eqnarray}
\|\nabla^2 u_t\|_{L^2}&\leq& C\|\mu\Delta
u_t+(\mu+\lambda)\nabla\divg
u_t\|_{L^2}\nonumber\\
&\leq&C\|\n u_{tt}-\rho_t \nabla f +\rho_t u_t+\rho_tu\cdot\nabla
u+\rho u_t\cdot\nabla u+\rho u\cdot\nabla
u_t+\nabla P_t\|_{L^2}\nonumber\\
&\leq&C\left(\|\rho^{1/2}
u_{tt}\|_{L^2}+\|\rho_t\|_{L^2}+\|\rho_t\|_{L^3}\|u_t\|_{L^6} +\|\rho_t\|_{L^3}\|u\|_{L^\infty}\|\nabla u\|_{L^6}\right.\nonumber\\
&&\left. +\|u_t\|_{L^6}\|\nabla u\|_{L^3}+\|u\|_{L^\infty}\|\nabla
u_t\|_{L^2}+\|\nabla P_t\|_{L^2}\right)\nonumber\\
&\leq& C\left(1+\|\nabla u_t\|_{L^2}+\|\rho^{1/2}
u_{tt}\|_{L^2}\right),\label{4.22}
\end{eqnarray}
from which and (\ref{4.21}) it follows that
$$
\int_0^T\sigma\|\nabla^2
u_t\|_{L^2}^2dt\leq C.
$$
This, together with (\ref{4.21}), finishes the proof of
(\ref{4.19}).\hfill$\square$

The next lemma is concerned with the $W^{1,q}$-estimate
($q\in(3,6)$) on the gradients of density and pressure, which in
particular indicates the H\"{o}lder continuity of $(\nabla
\rho,\nabla P)$.

\begin{lem}\label{lem4.5} Let $q\in(3,6)$ be as in Theorem \ref{thm1.1}. For any given $T>0$, it holds that
\begin{equation}
\sup_{0\leq t\leq T}\left(\|\nabla \rho\|_{W^{1,q}}+\|\nabla
P\|_{W^{1,q}}\right)+\int_0^T\left(\|u_t\|_{W^{1,q}}^{p_0}+\|\nabla^2
u\|_{W^{1,q}}^{p_0}\right)dt\leq C(T),\label{4.23}
\end{equation} where \be \la{a4.23}p_0\triangleq(9q-6)/(10q-12)\in
(1,4q/(5q-6)).\ee
\end{lem}
\pf Applying the differential operator $\na^2$ to both sides of
(\ref{4.13}), multiplying the resulting equations by $q|\na^2
P(\rho)|^{q-2}\na^2 P(\rho)$, and integrating it by parts over
$\R^3$, one deduces from (\ref{4.2}), (\ref{4.11}) and Lemma
\ref{lem2.1} that
\begin{eqnarray*}
\frac{d}{dt}\|\nabla^2 P\|_{L^q}^q&\leq& C\left(\|\nabla
u\|_{L^\infty}\|\nabla^2 P\|_{L^q}+\|\nabla P\|_{L^\infty}\|\nabla^2
u\|_{L^q}+\|\nabla^2u\|_{W^{1,q}}\right)\|\nabla^2
P\|_{L^q}^{q-1}\nonumber\\
&\leq&C\left(1+\|\nabla
u\|_{H^2}\right)\left(1+\|\nabla^2P\|_{L^q}^q\right)+C\|\nabla^2u\|_{W^{1,q}}\|\nabla^2P\|_{L^q}^{q-1}.
\end{eqnarray*}
The similar estimate also holds for $\|\nabla^2\rho\|_{L^q}$.
Therefore,
\begin{eqnarray}
\frac{d}{dt}\left(\|\nabla^2 \rho\|_{L^q}^q+\|\nabla^2
P\|_{L^q}^q\right)&\leq& C\left(1+\|\nabla
u\|_{H^2}\right)\left(1+\|\nabla^2\rho\|_{L^q}^q+\|\nabla^2P\|_{L^q}^q\right)\nonumber\\
&&+C\|\nabla^2u\|_{W^{1,q}}
\left(\|\nabla^2\rho\|_{L^q}^{q-1}+\|\nabla^2P\|_{L^q}^{q-1}\right).\label{4.24}
\end{eqnarray}

Applying the standard $W^{1,p}$-estimate to the  elliptic system  (\ref{1.2}) yields that
\begin{eqnarray}
\|\nabla^2 u\|_{W^{1,q}}&\leq& C\left(\|\rho u_t\|_{W^{1,q}}+\|\rho
u\cdot\nabla u\|_{W^{1,q}}+\|\nabla P\|_{W^{1,q}}+\|\rho\nabla
f\|_{W^{1,q}}\right)\nonumber\\
&\leq&C\left(1+\|u_t\|_{W^{1,q}}+\|\nabla\rho\|_{L^q}\|u_t\|_{L^\infty}+\|\nabla u\|_{W^{1,q}}\nonumber\right.\\
&&+\left.\|\nabla\rho\|_{L^\infty}\|\nabla u\|_{L^q}+\|\nabla
u\|_{L^\infty}\|\nabla u\|_{L^q}+\|\nabla
P\|_{W^{1,q}}\right)\nonumber\\
&\leq& C\left(1+\|\nabla
u\|_{H^2}+\|u_t\|_{W^{1,q}}+\|\nabla\rho\|_{W^{1,q}}+\|\nabla
P\|_{W^{1,q}}\right),\label{4.25}
\end{eqnarray}
where we have also used (\ref{4.2}), (\ref{4.11}) and Lemma
\ref{lem2.1}. Putting (\ref{4.25}) into (\ref{4.24}) gives
\begin{equation}
\frac{d}{dt}\left(\|\nabla^2 \rho\|_{L^q}^q+\|\nabla^2
P\|_{L^q}^q\right)\leq C\left(1+\|\nabla
u\|_{H^2}+\|u_t\|_{W^{1,q}}\right)\left(1+\|\nabla^2\rho\|_{L^q}^q+\|\nabla^2P\|_{L^q}^q\right).\label{4.26}
\end{equation}

Using (\ref{2.1}), (\ref{2.7}), (\ref{3.9}) and (\ref{4.10}), we
find that
\begin{eqnarray}
\int_0^T\|u_t\|_{L^q}^2dt&\leq& C\int_0^T\|u_t\|_{H^1}^2dt\leq
C\int_0^T\left(\|u_t\|_{L^2}^2+\|\nabla
u_t\|_{L^2}^2\right)dt\nonumber\\
&\leq& C\int_0^T\left(\|\rho^{1/2}u_t\|_{L^2}^2+\|\nabla
u_t\|_{L^2}^2\right)dt\leq C.\label{4.27}
\end{eqnarray}
Note that $ 4q/(5q-6)\in (1,4/3)$ for $q\in(3,6)$. So, for  $p_0 $ as in \eqref{a4.23}, we obtain by Lemma \ref{lem2.1}, H\"{o}lder
inequality and (\ref{4.19}) that
\begin{eqnarray}
&&\int_0^T\|\nabla u_t\|_{L^q}^{p_0}dt\leq C\int_0^T\|\nabla
u_t\|_{L^2}^{p_0(6-q)/(2q)}\|\nabla u_t\|_{L^6}^{p_0(3q-6)/(2q)}dt\nonumber\\
&&\quad\leq C\int_0^T\sigma^{-p_0/2}\left(\sigma\|\nabla
u_t\|_{L^2}^2\right)^{p_0(6-q)/(4q)}\left(\sigma\|\nabla
u_t\|_{H^1}^2\right)^{p_0(3q-6)/(4q)}dt\nonumber\\
&&\quad\leq C\left(\sup_{0\leq t\leq T}\sigma\|\nabla
u_t\|_{L^2}^2\right)^{p_0(6-q)/(4q)}\int_0^T\sigma^{-p_0/2}\left(\sigma\|\nabla
u_t\|_{H^1}^2\right)^{p_0(3q-6)/(4q)}dt\nonumber\\
&&\quad\leq C\left(\int_0^T
\sigma^{-2p_0q/(4q-p_0(3q-6))}dt\right)^{(4q-p_0(3q-6))/(4q)}\left(\int_0^T\sigma\|\nabla
u_t\|_{H^1}^2dt\right)^{p_0(3q-6)/(4q)}\nonumber\\
&&\quad\leq C,\label{4.28}
\end{eqnarray}
since $0<2p_0q/(4q-p_0(3q-6))<1$ and $0<p_0(3q-6)/(4q)<1$.

The combination of (\ref{4.27}) with (\ref{4.28}) shows that for
$q\in(3,6)$,
\begin{equation}
\int_0^T\|u_t\|_{W^{1,q}}^{p_0}dt\leq C .\label{4.29}
\end{equation}
Thus, by (\ref{4.11}), (\ref{4.29}) and Gronwall's inequality, one
sees from (\ref{4.26}) that
$$
\sup_{0\leq t\leq T}\left(\|\nabla\rho\|_{W^{1,q}}+\|\nabla
P\|_{W^{1,q}}\right)\leq C,
$$
which, combining with (\ref{4.11}), (\ref{4.25}), and (\ref{4.29}),
finishes the proof of Lemma \ref{lem4.5}.\hfill$\square$

Finally, we still need the following lemma, which implies that $u_t$
and $\nabla^2 u$ are H\"{o}lder continuous away from $t=0$.
\begin{lem}\label{lem4.6}For any given $T>0$, it holds that
\begin{equation}
\sup_{0\leq t\leq T}\sigma\left(\|\rho^{1/2}
u_{tt}\|_{L^2}+\|\nabla^2 u_t\|_{L^2}+\|\nabla^2
u\|_{W^{1,q}}\right)+\int_0^T\sigma^2\|\nabla u_{tt}\|_{L^2}^2dt\leq
C(T).\label{4.30}
\end{equation}
\end{lem}
\pf Differentiating (\ref{a4.20})   with
respect to $t$ gives
\begin{eqnarray*}
&&\rho u_{ttt}+\rho u\cdot\nabla u_{tt}-\mu\Delta
u_{tt}-(\mu+\lambda)\nabla\divg u_{tt}\\
&&\quad=2\divg (\rho u)u_{tt}+\divg(\rho u)_tu_t-2(\rho
u)_t\cdot\nabla u_t\\
&&\qquad-(\rho_{tt}u+2\rho_tu_t)\cdot\nabla u-\rho u_{tt}\cdot\nabla
u-\nabla P_{tt}+\rho_{tt}\nabla f,
\end{eqnarray*}
which, multiplied by $u_{tt}$ in $L^2$ and integrated by parts over
$\R^3$, yields
\begin{eqnarray}
&&\frac{1}{2}\frac{d}{dt}\int \rho|u_{tt}|^2dx+\int \left(\mu|\nabla
u_{tt}|^2+(\mu+\lambda)|\divg u_{tt}|^2\right)dx\nonumber\\
&&\quad=-4\int \rho u\cdot\nabla u_{tt}\cdot u_{tt}dx-\int (\rho
u)_t\cdot\left(\nabla(u_t\cdot
u_{tt})+2\nabla u_t\cdot u_{tt}\right)dx\nonumber\\
&&\qquad-\int \left(\rho_{tt} u+2\rho_tu_t\right)\cdot \nabla u\cdot
u_{tt}dx-\int \rho u_{tt}\cdot\nabla u\cdot
u_{tt}dx\nonumber\\
&&\qquad+\int P_{tt}\divg u_{tt}dx+\int \rho_{tt}\nabla f \cdot
u_{tt}dx\triangleq\sum_{i=1}^6J_i.\label{4.31}
\end{eqnarray}

The right-hand side of (\ref{4.31}) can be estimated term by term as
follows, using Lemma \ref{lem2.1}, Cauchy-Schwarz inequality and the
estimates obtained.
\begin{eqnarray*}
J_1&\leq&\|\rho^{1/2} u_{tt}\|_{L^2}\|\nabla u_{tt}\|_{L^2}\leq
\delta \|\nabla
u_{tt}\|_{L^2}^2+C(\delta)\|\rho^{1/2} u_{tt}\|_{L^2}^2,\\[2mm]
J_2&\leq&C\left(\|\rho
u_{t}\|_{L^3}+\|\rho_t\|_{L^3}\right)\left(\|\nabla
u_t\|_{L^2}\|u_{tt}\|_{L^6}+\|u_t\|_{L^6}\|\nabla
u_{tt}\|_{L^2}\right)\\
&\leq&C\left(1+\|\rho u_t\|_{L^2}^{1/2}\|\nabla
u_t\|_{L^2}^{1/2}\right)\|\nabla u_{t}\|_{L^2}\|\nabla
u_{tt}\|_{L^2}\\
&\leq& \delta\|\nabla u_{tt}\|_{L^2}^2+C(\delta)\left(1+\|\nabla
u_t\|_{L^2}^{3}\right),\\[2mm]
J_3&\leq& C\left(\|\rho_{tt}\|_{L^2}\|u\|_{L^\infty}\|\nabla
u\|_{L^3}+\|\rho_t\|_{L^3}\|u_t\|_{L^6}\|\nabla
u\|_{L^3}\right)\|u_{tt}\|_{L^6}\\
&\leq& \delta\|\nabla u_{tt}\|_{L^2}^2
+C(\delta)\left(\|\rho_{tt}\|_{L^2}^2+\|\nabla
u_t\|_{L^2}^2\right),\\[2mm]
J_4&\leq&C\|\rho^{1/2} u_{tt}\|_{L^2}\|\nabla
u\|_{L^3}\|u_{tt}\|_{L^6}\leq \delta\|\nabla
u_{tt}\|_{L^2}^2+C(\delta)\|\rho^{1/2} u_{tt}\|_{L^2}^2
\end{eqnarray*}
and
\begin{eqnarray*}
J_5+J_6&\leq&
C\left(\|P_{tt}\|_{L^2}+\|\rho_{tt}\|_{L^2}\right)\left(\|\nabla
u_{tt}\|_{L^2}+\|\nabla f\|_{L^3}\|u_{tt}\|_{L^6}\right)\\
&\leq&\delta\|\nabla
u_{tt}\|_{L^2}^2+C(\delta)\left(\|P_{tt}\|_{L^2}^2+\|\rho_{tt}\|_{L^2}^2\right).
\end{eqnarray*}

Putting the estimates of $J_1,\ldots,J_6$ into (\ref{4.31}),
multiplying it by $\sigma^2$ and integrating the resulting relation
over $(0,T)$,  we have by choosing $\delta>0$ small enough that
\begin{eqnarray}
&&\sup_{0\leq t\leq
T}\left(\sigma^2\|\rho^{1/2}u_{tt}\|_{L^2}^2\right)+\int_0^T\sigma^2\|\nabla
u_{tt}\|_{L^2}^2dt\nonumber\\
&&\quad\leq C+C\int_0^T\left(\sigma\|\rho^{1/2}
u_{tt}\|_{L^2}^2+\sigma^2\left(\|P_{tt}\|_{L^2}^2+\|\rho_{tt}\|_{L^2}^2+\|\nabla
u_t\|_{L^2}^3\right)\right)dt\nonumber\\
&&\quad\leq C+C\sup_{0\leq t\leq T}\left(\sigma^{1/2}\|\nabla
u_t\|_{L^2}\right)\int_0^T\|\nabla u_t\|_{L^2}^2dt\leq
C,\label{4.32}
\end{eqnarray}
where we have used (\ref{4.10}), (\ref{4.12}) and (\ref{4.19}).

As a result of (\ref{4.32}), (\ref{4.19}) and (\ref{4.22}), we also
see that
\begin{equation}
\sigma\|\nabla^2 u_t\|_{L^2} \leq C\sigma\left(1+\|\nabla
u_t\|_{L^2}+\|\rho^{1/2} u_{tt}\|_{L^2}\right)\leq C,\quad\forall\;
t\in[0,T].\label{4.33}
\end{equation}
and thus, it follows from (\ref{4.25}), (\ref{4.19}) and
(\ref{4.23}) that for any $t\in[0,T]$,
\begin{eqnarray}
\sigma\|\nabla^2 u\|_{W^{1,q}}&\leq& C\sigma\left(1+\|\nabla
u\|_{H^2}+\|u_t\|_{W^{1,q}}+\|\nabla\rho\|_{W^{1,q}}+\|\nabla
P\|_{W^{1,q}}\right)\nonumber\\
&\leq&C\left(1+\sigma\|u_t\|_{H^2}\right)\leq C,\label{4.34}
\end{eqnarray}
where we have used (\ref{2.7}), (\ref{3.9}), (\ref{4.10}) and
(\ref{4.19}) to get that $\sigma\|u_t\|_{L^2}\leq C$. Lemma
\ref{lem4.6} now readily follows from
(\ref{4.32})--(\ref{4.34}).\hfill$\square$

\section{Proofs of Theorems \ref{thm1.1} and \ref{thm1.2}}
With all the a priori estimates at hand, we are now ready to prove
our main results. First of all, we prove the global well-posedness
of classical solutions to the problem (\ref{1.1})--(\ref{1.5})
provided the initial density is strictly away from vacuum and the
initial energy is small.
\begin{pro}\label{pro5.1} Let $C_0$ be the initial energy defined by (\ref{1.10}). For given numbers $M>0$ (not necessarily small)
and $\tilde\rho>0$, assume that $\rho_0,u_0$ and $f$ satisfy
\begin{equation}
\left\{
\begin{array}{l}
f\in H^3,\quad (\rho_0-\rho_\infty,u_0)\in H^3,\\[3mm]
0<\inf\limits_{x\in\R^3}\rho_0(x)\leq\sup\limits_{x\in\R^3}\rho_0(x)<\tilde\rho,\quad
\|\nabla u_0\|_{L^2}\leq M,\\[3mm]
C_0\leq \tilde\varepsilon,
\end{array}\right.\label{5.1}
\end{equation}
where $\tilde\varepsilon>0$ is the same one as in Proposition
\ref{pro3.1}. Then for any $0<T<\infty$, there exists a unique
classical solution $(\rho, u)$ of (\ref{1.1})--(\ref{1.5}) on
$\R^3\times(0,T]$ satisfying (\ref{2.9}), (\ref{2.10}) with $T_0$
replaced by any $T>0$. Moreover, all the uniform-in-time estimates
in Sect. 3 hold for $(\rho,u)$.
\end{pro}
\pf The standard local existence result (i.e. Lemma \ref{lem2.5})
shows that the Cauchy problem (\ref{1.1})--(\ref{1.5}) has a unique
local classical solution $(\rho, u)$, defined up to a positive $T_0$
which may depend on $\inf\limits_{x\in\R^3}\rho_0(x)$ and satisfying
(\ref{2.9}), (\ref{2.10}).

In view of (\ref{5.1}) and the definitions of $\Phi_i(T)$
$(i=1,2,3)$, we know that
$$
\Phi_1(0)=\Phi_2(0)=0,\quad\Phi_3(0)\leq M\quad{\rm and}\quad
\rho_0\leq \tilde\rho.
$$
Thus, there exists a positive $T_1\in(0,T_0]$ such that (\ref{3.4})
holds for $T=T_1$.

Set
\begin{equation}
T_*=\sup\left\{T\;\left|\;(\ref{3.4})\;\; {\rm
holds}\right.\right\}.\label{5.2}
\end{equation}
It is clear that $T_*\geq T_1>0$.

We claim that
\begin{equation}
T_*=\infty.\label{5.3}
\end{equation}
Otherwise, $T_*<\infty$. Then it follows from Proposition
\ref{pro3.1} that (\ref{3.5}) holds for any $0<T<T_*$, and
furthermore, the estimates in Lemmas \ref{lem4.1}--\ref{lem4.6} are
also valid for all $0<T<T_*$.

In the next, for the sake of simplicity, we denote by $\tilde C$ the
positive constant which may depend on the lower bound of the initial
density (i.e. $\inf\rho_0(x)$) and $T_*$. We shall prove that there
exists a positive constant $\tilde C>0$ such that for any
$T\in(0,T_*)$,
\begin{equation}
\sup_{0\leq t\leq T}\|(\rho-\rho_\infty,u)\|_{H^3}\leq \tilde
C.\label{5.4}
\end{equation}
This in particular implies that
$$
\|(\rho-\rho_\infty,u)(\cdot,T_*)\|_{H^3}<\infty,\quad
\inf_{x\in\R^3}\rho(x,T_*)>0.
$$
So, Lemma \ref{lem2.5}, together with (\ref{3.5}), yields that there
exists some $T_{**}>T_*$ such that (\ref{3.4}) holds for $T=T_{**}$,
which contradicts (\ref{5.2}). Hence, (\ref{5.3}) holds. Thus,
$(\rho,u)$ is in fact the unique classical solution of
(\ref{1.1})--(\ref{1.5}) on $\R^3\times(0,T]$ for all
$0<T<T_*=\infty$.

We are now in a position of proving (\ref{5.4}). To do so, we first
note that under the conditions of (\ref{5.1})$_1$ and
(\ref{5.1})$_2$, it holds that (keeping in mind that $\rho_0>0$)
\begin{equation}
u_t(\cdot,0)=\nabla f-u_0\cdot\nabla u_0+\rho_0^{-1}\left(\mu\Delta
u_0+(\mu+\lambda)\nabla\divg u_0-P(\rho_0)\right)\in H^1.\label{5.5}
\end{equation}
Thus, similar to the proof of (\ref{4.21}), by (\ref{5.5}) one
easily deduces that
\begin{equation}
\sup_{0\leq t\leq T}\left(\|\nabla u\|_{H^2}^2+\|\nabla
u_t\|_{L^2}^2\right)+\int_0^T\|\rho^{1/2}u_{tt}\|_{L^2}^2dt\leq
\tilde C.\label{5.6}
\end{equation}
Similar to the derivation of (\ref{4.22}), one gets  by using
(\ref{5.6}) that
\begin{eqnarray}
\int_0^T\|\nabla^2 u_t\|_{L^2}^2dt &\leq&\tilde C\int_0^T
\left(1+\|\nabla u_t\|_{L^2}^2+\|\rho^{1/2}
u_{tt}\|_{L^2}^2\right)dt\nonumber\\
&\leq& \tilde C+\tilde C\int_0^T\|\rho^{1/2}
u_{tt}\|_{L^2}^2dt\leq\tilde C.\label{5.7}
\end{eqnarray}
Moreover, using the standard $H^2$-theory of elliptic system, one
deduces from Lemmas \ref{lem2.1}, \ref{lem4.1}--\ref{lem4.3} and
(\ref{5.6}) that
\begin{eqnarray}
\|\nabla^2u\|_{H^2}&\leq& \tilde C\|\mu\Delta
u+(\mu+\lambda)\nabla\divg u\|_{H^2}\nonumber\\
&\leq& \tilde C\|\rho u_t+\rho u\cdot\nabla u+\nabla P-\rho\nabla
f\|_{H^2}\nonumber\\
&\leq&\tilde C\left(\|\rho
u_t\|_{L^2}+\|\nabla\rho\|_{H^1}\|u_t\|_{L^\infty}+\|\nabla\rho\|_{L^3}\|
\nabla u_t\|_{L^6}+\|\nabla
u_t\|_{H^1}\right.\nonumber\\
&&+\|\rho u\|_{H^2}\|\nabla u\|_{H^2}\left.+\|\nabla
P\|_{H^2}+\|\nabla f\|_{H^2}+\|\nabla
\rho\|_{H^1}\|\nabla f\|_{H^2}\right)\nonumber\\
&\leq&\tilde
C\left(1+\|\nabla^3\rho\|_{L^2}+\|\nabla^3P\|_{L^2}+\|\nabla
u_t\|_{H^1}\right).\label{5.8}
\end{eqnarray}

Applying the differential operator $\nabla^3$ to both sides of
(\ref{1.1}) and (\ref{4.13}), and multiplying the resulting
equations by $\nabla^3\rho$ and $\nabla^3P(\rho)$ in $L^2$,
respectively, one easily obtains after integrating by parts and
using Lemmas \ref{lem2.1}, \ref{lem4.1}--\ref{lem4.3}, (\ref{5.6})
that
\begin{equation}
\frac{d}{dt}\left(\|\nabla^3 \rho\|_{L^2}^2+\|\nabla^3
P\|_{L^2}^2\right)\leq\tilde
C\left(1+\|\nabla^2u\|_{H^2}^2+\|\nabla^3\rho\|_{L^2}^2+\|\nabla^3
P\|_{L^2}^2\right).\label{5.9}
\end{equation}
Combining (\ref{5.7}), (\ref{5.8}) and (\ref{5.9}), we conclude from
Gronwall's inequality that
\begin{equation}
\sup_{0\leq t\leq T}\left(\|\nabla^3 \rho\|_{L^2}^2+\|\nabla^3
P\|_{L^2}^2\right)\leq \tilde C.\label{5.10}
\end{equation}
So, Lemmas \ref{lem4.1}--\ref{lem4.3}, together with (\ref{5.6}) and
(\ref{5.10}), lead to (\ref{5.4}). The proof of Proposition
\ref{pro5.1} is therefore complete.\hfill$\square$

\vskip 2mm

{\it Proof of Theorem \ref{thm1.1}.} Let $(\rho_0,u_0)$ be the
initial data satisfying the conditions (\ref{1.12})--(\ref{1.14}) in
Theorem \ref{thm1.1}. Assume that the initial energy $C_0$ satisfies
\begin{equation}
C_0\leq \varepsilon\triangleq\tilde\varepsilon/2,\label{5.11}
\end{equation}
where $\tilde\varepsilon>0$ is the same one as in Proposition
\ref{pro3.1}.

To apply Proposition \ref{pro5.1}, we define the smooth approximate
data as follows:
\begin{equation}
\rho_0^{\delta,\eta}=\frac{j_\delta\ast\rho_0+\eta}{1+\eta},\quad
g^{\delta,\eta}=j_\delta\ast g,\quad f^{\delta,\eta}=j_\delta\ast
f,\label{5.12}
\end{equation}
where $0<\delta,\eta<1$, $j_\delta(x)$ is the standard mollifier
with width $\delta$, and `` $\ast$ " denotes the usual convolution
operator. As that in \cite{CK2003}, let $u_0^{\delta,\eta}$ be the
unique solution to the elliptic problem:
\begin{equation}
\left\{
\begin{array}{l} -\mu \Delta
u_0^{\delta,\eta}-(\mu+\lambda)\nabla\divg u_0^{\delta,\eta}=-\nabla
P(\rho_0^{\delta,\eta})+\left(\rho_0^{\delta,\eta}\right)^{1/2}g^{\delta,\eta},\\[2mm]
u_0^{\delta,\eta}(x)\to0\quad{\rm as}\quad |x|\to\infty.
\end{array}\right.\label{5.13}
\end{equation}
Also let $\rho_s^{\delta,\eta}\triangleq\rho_s^{\delta,\eta}(x)$ be
the unique solution of (\ref{1.6}) with $f$ being replaced by
$f^{\delta,\eta}$. Obviously,  $\rho_s^{\delta,\eta}$ satisfies
(\ref{1.8}), (\ref{1.9}). Moreover, it is easy to check that
\begin{equation}
\left\{
\begin{array}{l}
0<\frac{\eta}{1+\eta}\leq \rho_0^{\delta,\eta}\leq
\tilde\rho<\infty,\\[2mm]
(\rho_0^{\delta,\eta}-\rho_\infty,g^{\delta,\eta},u_0^{\delta,\eta},
f^{\delta,\eta},\rho_s^{\delta,\eta})\in H^\infty,\\[2mm]
\lim\limits_{\delta,\eta\to0}\left(\|\rho_0^{\delta,\eta}-\rho_0\|_{H^2{\cap}W^{2,q}}+\|u_0^{\delta,\eta}-u_0\|_{H^2}\right)=0,\\[2mm]
\lim\limits_{\delta,\eta\to0}\left(\|g^{\delta,\eta}-g\|_{L^2}+\|(f^{\delta,\eta}-f,\rho_s^{\delta,\eta}-\rho_s)\|_{H^2{\cap}W^{2,q}}\right)=0.
\end{array}
\right.\label{5.14}
\end{equation}
The initial norm of the mollified data
$(\rho_0^{\delta,\eta},u_0^{\delta,\eta})$ now reads
$$
C_0^{\delta,\eta}\triangleq \int
\left(G(\rho_0^{\delta,\eta})+\frac{1}{2}\rho_0^{\delta,\eta}|u_0^{\delta,\eta}|^2\right)dx,
$$
where $G(\cdot)$ is the function defined in (\ref{1.11}) with
$\rho_s$ replaced by $\rho_s^{\delta,\eta}$.

By (\ref{5.14}), we have
$$
\lim_{\eta\to0}\lim_{\delta\to0}C_0^{\delta,\eta}=C_0.
$$
So, there exist positive constants $\eta_0\in(0,1)$ and
$\delta_0(\eta)$ such that
\begin{equation}
C_0^{\delta,\eta}\leq C_0+\tilde\varepsilon/2\leq
\tilde\varepsilon,\label{5.15}
\end{equation}
provided that
\begin{equation}
0<\eta<\eta_0\quad{\rm and}\quad0<\delta<\delta_0(\eta).\label{5.16}
\end{equation}

Let $\delta,\eta$ satisfy (\ref{5.16}). Then it follows from
(\ref{5.15}) and Proposition \ref{pro5.1} that there exists a unique
smooth solution $(\rho^{\delta,\eta},u^{\delta,\eta})$ of
(\ref{1.1}), (\ref{1.2}) with the mollified data
$\rho_0^{\delta,\eta},u_0^{\delta,\eta},f^{\delta,\eta}$ and
$g^{\delta,\eta}$ on $\R^3\times (0,T]$ for all $T>0$. Moreover,
Proposition \ref{pro3.1} and Lemmas \ref{lem4.1}--\ref{lem4.6},
independent of $\delta$ and $\eta$, hold for
$(\rho^{\delta,\eta},u^{\delta,\eta})$.

Now, passing to the limit first $\delta\to0$, then $\eta\to0$, we
deduce from standard arguments that there exists a unique solution
$(\rho,u)$ to the origin problem (\ref{1.1})--(\ref{1.5}) on
$\R^3\times(0,T]$ for all $T>0$ satisfying
\begin{equation}
\left\{
\begin{array}{lll}
0\leq\rho(x,t)\leq 2\tilde\rho\quad{\rm for\;\;all}\quad
x\in\R^3,t\geq0,\\[2mm]
\left(\rho-\rho_\infty,P(\rho)-P(\rho_\infty)\right)\in L^\infty(0,T;H^2\cap W^{2,q}),\\[2mm]
u\in L^\infty(0,T;H^2)\cap L^2(0,T;H^3)\cap L^\infty(\tau,T; H^3\cap W^{3,q}),\\[2mm]
u_t\in L^\infty(\tau,T;H^2)\cap H^1(\tau,T;H^1),
\end{array}\right.\label{5.17}
\end{equation}
for any $0<\tau<T<\infty$. Note that, the uniqueness of $(\rho,u)$
satisfying (\ref{5.17}) can be proved in a standard way as that in
\cite{CK2003}.

To complete the proof of the first part of Theorem \ref{thm1.1}, we
still need to show that
\begin{equation}
\left(\rho-\rho_\infty,P(\rho)-P(\rho_\infty)\right)\in
C([0,T];H^2\cap W^{2,q}),\quad u\in C([0,T];H^2).\label{5.18}
\end{equation}
The proof of \eqref{5.18} is
similar to that in \cite{HuLi}, and we sketch it here for
completeness.
First, by virtue of (\ref{4.2}), (\ref{4.10})--(\ref{4.12}) and
(\ref{4.23}), it is easy to get that
\begin{equation}
\left\{
\begin{array}{l}
\left(\rho-\rho_\infty,P(\rho)-P(\rho_\infty)\right)\in
C([0,T];H^1\cap W^{1,q}),\\[2mm]
\left(\rho-\rho_\infty,P(\rho)-P(\rho_\infty)\right)\in C([0,T];H^2\cap W^{2,q}-{\rm
weak}),\\[2mm]
u\in C([0,T];H^1\cap W^{1,q}).
\end{array}\right.\label{5.19}
\end{equation}
Denote $D_{ij}\triangleq\partial_{ij}^2$ with $i,j=1,2,3$. Note that
it holds in $\mathcal{D}'(\R^3\times(0,T))$ that
$$
\pa_t D_{ij}\rho+\divg(uD_{ij}\rho)=-\divg(\rho
D_{ij}u)-\divg(\pa_i\rho\cdot\pa_ju+\pa_j\rho\cdot\pa_iu).
$$
Let $j_\nu(x)$ be the standard mollifying kernel with width $\nu$,
and set $\rho^\nu\triangleq\rho\ast j_\nu$. Then,
\begin{equation}
\pa_t D_{ij}\rho^\nu+\divg(uD_{ij}\rho^\nu)=-\divg(\rho D_{ij}u)\ast
j_\nu-\divg(\pa_i\rho\cdot\pa_ju+\pa_j\rho\cdot\pa_iu)\ast j_\nu
+R_\nu,\label{5.20}
\end{equation}
where $R_\nu\triangleq\divg(uD_{ij}\rho^\nu)-\divg(uD_{ij}\rho)\ast
j_\nu$ satisfies (cf. \cite[Lemma 2.3]{Li1996})
\begin{equation}
\int_0^T\|R_\nu\|_{L^2\cap L^q}^{p_0}dt\leq
C\int_0^T\|u\|_{W^{1,\infty}}^{p_0}\|D_{ij}\rho\|_{L^2\cap
L^q}^{p_0}dt\leq C,\label{5.21}
\end{equation} with $p_0>1$   as in (\ref{a4.23}).

Multiplying (\ref{5.20}) by $q|D_{ij}\rho^\nu|^{q-2}D_{ij}\rho^\nu$
and integrating by parts over $\R^3$, we obtain
\begin{eqnarray*}
\frac{d}{dt}
\|D_{ij}\rho^\nu\|_{L^q}^q&=&-(q-1)\int|D_{ij}\rho^\nu|^q\divg u
dx-q\int\divg(\rho D_{ij}u)\ast j_\nu
|D_{ij}\rho^\nu|^{q-2}D_{ij}\rho^\nu dx\nonumber\\
&&-q\int\left(\divg(\pa_i\rho\cdot\pa_ju+\pa_j\rho\cdot\pa_iu)\ast
j_\nu\right) |D_{ij}\rho^\nu|^{q-2}D_{ij}\rho^\nu
dx\nonumber\\
&&+q\int R_\nu |D_{ij}\rho^\nu|^{q-2}D_{ij}\rho^\nu dx,
\end{eqnarray*}
which, combining with (\ref{4.11}), (\ref{4.23}) and (\ref{5.21}),
yields
\begin{eqnarray*}
&&\sup_{0\leq t\leq T} \|\nabla^2\rho^\nu\|_{L^q}^q+
\int_0^T\left|\frac{d}{dt} \|\nabla^2\rho^\nu\|_{L^q}^q\right|^{p_0}dt\\
&&\quad\leq C+C\int_0^T\left(\|\nabla
u\|_{W^{2,q}}^{p_0}+\|R_\nu\|_{L^q}^{p_0}\right)dt\leq C.
\end{eqnarray*}
This, together with Ascoli-Arzela theorem, gives
$$
\|\nabla^2\rho^\nu(\cdot,t)\|_{L^q}\to
\|\nabla^2\rho(\cdot,t)\|_{L^q}\quad {\rm in}\quad
C([0,T]),\quad{\rm as}\quad \nu\to0,
$$
which particularly implies
\begin{equation}
\|\nabla^2\rho(\cdot,t)\|_{L^q}\in C([0,T]).\label{5.22}
\end{equation}
Similarly, one can also obtain that
\begin{equation}
\|\nabla^2\rho(\cdot,t)\|_{L^2}\in C([0,T]).\label{5.23}
\end{equation}
So, it readily follows from (\ref{5.22}), (\ref{5.23}) and
(\ref{5.19})$_2$ that
\begin{equation}
\nabla^2\rho \in C([0,T];L^2\cap L^q).\label{5.24}
\end{equation}
In the exactly same way, we also have
\begin{equation}
\nabla^2P(\rho) \in C([0,T];L^2\cap L^q).\label{5.25}
\end{equation}

It follows from (\ref{1.2}) that
$$
(\rho u_t)_t=\left(\mu\Delta u_t+(\mu+\lambda)\nabla\divg
u_t\right)+\rho_t\nabla f-\nabla P_t-(\rho u\cdot\nabla u)_t.
$$
Thus, it is easily seen from (\ref{4.2}), (\ref{4.10}) and
(\ref{4.12}) that
$$
\|(\rho u_t)_t\|_{H^{-1}}\leq C\left(1+\|\nabla u_t\|_{L^2}\right),
$$
and hence,
$$
(\rho u_t)_t\in L^2(0,T;{H^{-1}}).
$$
This, combining with the fact that $\rho u_t\in L^2(0,T;H^1)$ due to
(\ref{4.10}) and (\ref{4.11}), leads to
\begin{equation}
\rho u_t\in C([0,T]; L^2).\label{5.26}
\end{equation}
Thus, by virtue of (\ref{5.19}), (\ref{5.24}) and (\ref{5.25}), one
easily gets
\begin{equation}
u\in C([0,T];H^2),\label{5.27}
\end{equation}
since it holds that
$$
\mu\Delta u+(\mu+\lambda)\nabla\divg u=\rho u_t+\rho u\cdot\nabla
u+\nabla P(\rho)-\rho\nabla f.
$$
Combining (\ref{5.19}), (\ref{5.24}), (\ref{5.25}), and
(\ref{5.27}), we finish  the proof of (\ref{5.18}).

To complete the proof of Theorem \ref{thm1.1}, it remains to  prove
(\ref{1.18}). To do this, we first deduce
in a manner similar to the derivation of (\ref{3.35}) that
\begin{eqnarray}
&&\int_1^\infty\left|\frac{d}{dt}\|\rho-\rho_s\|_{L^4}^4\right|dt\nonumber\\
&&\quad\leq C\int_1^\infty\left(\|\nabla
u\|_{L^2}^2+\|\rho-\rho_s\|_{L^4}^4+\|F\|_{L^4}^4\right)dt\leq
C,\label{5.28}
\end{eqnarray}
where we have used (\ref{3.7}) and (\ref{3.34}). Thus, it follows
from (\ref{3.34}) and (\ref{5.28}) that
$$
\|\rho-\rho_s\|_{L^4}\to0\quad{\rm as}\quad t\to\infty.
$$
This, together with (\ref{3.9}), Proposition \ref{pro1.1} and the
interpolation inequality, immediately gives
\begin{equation}
\|\rho-\rho_\infty\|_{L^p}\to0\quad{\rm as}\quad
t\to\infty,\quad\forall\;p\in(2,\infty).\label{5.29}
\end{equation}

To study the large-time behavior of the velocity, we set
$$
M(t)\triangleq\frac{\mu}{2}\|\nabla
u\|_{L^2}^2+\frac{\mu+\lambda}{2}\|\divg u\|_{L^2}^2.
$$
Then, multiplying (\ref{3.12}) by $\dot u$ in $L^2$ and integrating
by parts, we obtain
\begin{equation}
|M'(t)|\leq C\int \rho|\dot u|^2dx+C\|\nabla
u\|_{L^3}^3+C\|\nabla\dot u\|_{L^2}.\label{5.30}
\end{equation}
Here, we have used (\ref{3.14}), (\ref{3.15}) and the following
estimates:
$$
\left|\int\dot u\cdot\nabla\left(P(\rho)-P(\rho_s)\right)
dx\right|\leq C\|\rho-\rho_s\|_{L^2}\|\nabla\dot u\|_{L^2}\leq
C\|\nabla\dot u\|_{L^2}
$$
and
$$
\left|\int(\rho-\rho_s)\dot u\cdot\nabla f dx\right|\leq
C\|\rho-\rho_s\|_{L^2}\|\dot u\|_{L^6}\|\nabla f\|_{L^3}\leq
C\|\nabla\dot u\|_{L^2}.
$$
Thus, it follows from (\ref{3.34}), (\ref{3.39}), (\ref{3.41}) and
(\ref{5.30}) that
\begin{eqnarray}
\int_1^\infty|M'(t)|^2dt&\leq& C\int_1^\infty\left(\|\rho^{1/2}\dot
u\|_{L^2}^4+\|\nabla u\|_{L^2}^2\|\nabla u\|_{L^4}^{4}+\|\nabla\dot
u\|_{L^2}^2\right)dt\nonumber\\
&\leq& C\int_1^\infty\left(\|\rho^{1/2}\dot u\|_{L^2}^2+\|\nabla
u\|_{L^4}^{4}+\|\nabla\dot u\|_{L^2}^2\right)dt\leq C.\label{5.31}
\end{eqnarray}

Due to (\ref{3.40}) and (\ref{3.7}), one has
$$
\int_1^\infty|M(t)|^2dt\leq C\sup_{t\geq1}\|\nabla
u\|_{L^2}^2\int_1^\infty\|\nabla u\|_{L^2}^2dt\leq C,
$$
from which and (\ref{5.31}), we know that
\begin{equation}
\|\nabla u(t)\|_{L^2}\to0\quad{\rm as}\quad t\to\infty.\label{5.32}
\end{equation}
As a result, we also have
$$
\int\rho^{1/2}|u|^4dx\leq \left(\int\rho
|u|^2dx\right)^{1/2}\|u\|_{L^6}\leq C\|\nabla u\|_{L^2}\to0\quad{\rm
as}\quad t\to\infty.
$$
This, together with (\ref{5.29}) and (\ref{5.32}), proves
(\ref{1.18}). Therefore, the proof of Theorem \ref{thm1.1} is
complete.\hfill$\square$

\vskip 2mm

{\it Proof of Theorem \ref{thm1.2}.} Define the approximate data
$(\rho_0^{\delta,\eta},u_0^{\delta,\eta},f^{\delta,\eta})$  as
follows:
$$
\rho_0^{\delta,\eta}=\frac{j_\delta\ast\rho_0+\eta}{1+\eta},\quad
u_0^{\delta,\eta}=j_\delta\ast u_0,\quad
f^{\delta,\eta}=j_\delta\ast f.
$$
Then, we can apply Proposition \ref{pro5.1} to obtain a global
smooth solution $(\rho^{\delta,\eta},u^{\delta,\eta})$ of the Cauchy
problem (\ref{1.1}), (\ref{1.2}) with the mollified data
$\rho_0^{\delta,\eta},u_0^{\delta,\eta}$ and $f^{\delta,\eta}$,
which satisfies the uniform bounds in Proposition \ref{pro3.1}. Now,
the remaining arguments to obtain the global weak solution and its
asymptotic behavior are almost the same as that of  \cite{Ho3}. The
proof of Theorem \ref{thm1.2} is thus finished.\hfill$\square$

\end{document}